\documentclass[11pt]{article}

\usepackage[utf8]{inputenc}
\usepackage{lipsum}
\usepackage[dvipsnames]{xcolor}
\usepackage{url}
\usepackage{bm}
\usepackage{mathpazo}
\usepackage[T1]{fontenc}
\usepackage[utf8]{inputenc}
\usepackage[english]{babel}
\usepackage[active]{srcltx}
\usepackage{colorspace}
\usepackage{tikz}
\usepackage{sansmath}
\usetikzlibrary{shadings,intersections,patterns}
\usepackage{pgfplots}

\pgfplotsset{compat=1.10}
\usepgfplotslibrary{fillbetween}
\usetikzlibrary{patterns}
\usetikzlibrary{intersections, calc, angles}
\usepackage{tkz-euclide}
\pgfdeclarelayer{pre main}
\usepackage{environ}
\makeatletter
\newsavebox{\measure@tikzpicture}
\NewEnviron{scaletikzpicturetowidth}[1]{%
  \def\tikz@width{#1}%
  \begin{lrbox}{\measure@tikzpicture}%
  \BODY
  \end{lrbox}%
  \pgfmathparse{#1/\wd\measure@tikzpicture}%
  \BODY
}
\makeatother
\usetikzlibrary{arrows,automata}
\usetikzlibrary{positioning}
\usetikzlibrary{calc,through}
\usepackage{etoolbox}%
\apptocmd{\thebibliography}{\fontsize{11}{15}\selectfont}{}{}%

\tikzset{
    state/.style={
           rectangle,
           rounded corners,
           draw=black, very thick,
           minimum height=2em,
           inner sep=2pt,
           text centered,
           },
}

\usepackage{mathrsfs}

\usepackage[left=2.5cm, right=2.5cm, top=3cm, bottom=3cm]{geometry}
\usepackage{amsmath}
\usepackage{amssymb}
\usepackage{amsthm}
\usepackage{graphicx}
\usepackage{subcaption}
\usepackage{caption}
\usepackage{xcolor}
\usepackage{longtable}
\usepackage{verbatimbox}
\usepackage{xfrac}

\usepackage[shortlabels]{enumitem}

\usepackage{hyperref} 
\usepackage[capitalise, nameinlink, noabbrev]{cleveref} 
\hypersetup{colorlinks=true, 
	linkcolor=blue,
	citecolor=red,
}

\usepackage[title]{appendix}
\usepackage{bm}
\usepackage{empheq}
\usepackage{placeins}
\theoremstyle{plain}
\newtheorem{theorem}{Theorem}[section]

\newtheorem{remark}[theorem]{Remark}

\theoremstyle{definition}
\theoremstyle{remark}
\numberwithin{equation}{section}

\theoremstyle{plain}
\newtheorem*{theorem*}{Theorem}
\newtheorem*{corollary*}{Corollary}

\theoremstyle{definition}

\newtheorem*{notation*}{Notation}

\numberwithin{equation}{section}
\numberwithin{figure}{section}

\newcommand{\R}{{\mathbb R}}

\newcommand{\N}{{\mathbb N}}
\newcommand{\Z}{{\mathbb Z}}
\newcommand{\vect}[1]{\boldsymbol{#1}}

\newcommand{\e}{\varepsilon}
\newcommand{\res}{\mathop{\hbox{\vrule height 7pt width .5pt depth 0pt
\vrule height .5pt width 6pt depth 0pt}}\nolimits}

\newcommand{\boundellipse}[3]
{(#1) ellipse (#2 and #3)
}

\usepackage[small, sc]{titlesec}

\usepackage[
backend=biber,
style=numeric,
sorting=nty,
maxbibnames=50,
giveninits=true
]{biblatex}
\usepackage{csquotes}
\renewbibmacro{in:}{%
 \ifentrytype{article}{}{}}
\DeclareFieldFormat[article]{volume}{\textbf{#1}}
\DeclareFieldFormat[article]{title}{\textit{#1}}
\DeclareFieldFormat[article]{journaltitle}{#1}
\DeclareFieldFormat[article]{number}{\textbf{#1}}
\renewbibmacro*{volume+number+eid}{%
  \printfield{volume}%
\setunit*{\textbf{\addcolon}}
  \printfield{number}\setunit{\addcomma\space}%
  \printfield{eid}}

\usepackage{xcolor}
\hypersetup{
colorlinks,
linkcolor={red!50!black},
citecolor={blue!50!black},
urlcolor={blue!80!black}
}

\begin{document}
\title{\textsc{Soap films: from the Plateau problem\\ to deformable boundaries}}

\author{\textsc{Giulia Bevilacqua}$^1$\thanks{\href{mailto:giulia.bevilacqua@dm.unipi.it}{\texttt{giulia.bevilacqua@dm.unipi.it}}}\,\,\,$-$\,\, \textsc{Luca Lussardi}$^2$\thanks{\href{mailto:luca.lussardi@polito.it}{
\texttt{luca.lussardi@polito.it}}} \,\,\,$-$\,\, \textsc{Alfredo Marzocchi}$^3$\thanks{\href{mailto:alfredo.marzocchi@unicatt.it}{
\texttt{alfredo.marzocchi@unicatt.it}}}\bigskip\\
\normalsize$^1$ Dipartimento di Matematica, Università di Pisa, Largo Bruno Pontecorvo 5, I–56127 Pisa, Italy\\
\normalsize$^2$ DISMA ``G.L.\,Lagrange'', Politecnico di Torino, c.so Duca degli Abruzzi 24, I-10129 Torino, Italy\\
\normalsize$^3$ Dipartimento di Matematica e Fisica ``N.\,Tartaglia'', Università Cattolica del Sacro Cuore,\\
\normalsize via della Garzetta 48, I-25133 Brescia, Italy\\
}

\date{}

\maketitle

\begin{abstract}
\noindent A review on the classical Plateau problem is presented. Then, the state of the art about the Kirchhoff-Plateau problem is illustrated as well as some possible future directions of research. 
\end{abstract}

\bigskip



\bigskip
\bigskip

\maketitle

\section{Introduction}

Soap films arise as equilibrium interfaces between two fluids. Young and Laplace, at the beginning of the 19th century, gave an expression for the pressure difference $p_c$ (the {\it capillary pressure}) over an interface between two fluids. Precisely, the relation, which is now called {\it Young-Laplace equation}, is given by
\[
p_c=\sigma \left(\frac{1}{r_1}+\frac{1}{r_2}\right),
\]
where $\sigma>0$ is a constant called {\it surface tension} and $r_1,r_2$ are the principal radii of curvature of the interface $S$ (see, for instance, \cite{SS} for a derivation of the Young-Laplace equation). The surface tension measures the amount of energy one needs to extend the surface $S$ by one unit area. Taking into account the definition of mean curvature of $S$, denoted by $H$, we can rewrite the Young-Laplace equation as $p_c=\sigma H$. As a consequence, the interface is in equilibrium if and only if $H$ is constant. From the physical point of view, two different kind of configurations are essentially possible. In one case, the interface $S$ forms a closed surface, that is a compact surface without boundary; then, $S$ must be a sphere, and this explains why soap bubbles are round. In the other case, the interface $S$ is a surface with boundary. In this case, the Young-Laplace equation becomes $H=0$. The best physical model for these kind of surfaces is represented by soap films: putting a rigid wire in a soap solution and extracting it, a thin soap film will remains attached to the wire. In the middle of the 19th century, the Belgian physicist Plateau devised many experiments putting rigid wires in a soap solution in order to understand the possible singular configurations of soap films. For this reason, still today we use the terminology {\it Plateau problem} to deal with the problem of finding the shape of soap films with some prescribed boundaries. 

From a mathematical point of view, soap films turn out to be stable {\it minimal surfaces}. The connection between soap films and minimal surfaces dates back to Gauss who worked, in the 19th century, on capillarity problems. Actually, in the middle of the 18th century, Lagrange was the first who investigated minimal surfaces as critical points of the area functional. Indeed, at least in the smooth case, the {\it minimal surface equation} $H=0$ is the Euler-Lagrange equation of the area functional. This suggests an interesting change of point of view: instead of solving directly the partial differential equation $H=0$, one might look at minimizers of the area functional. This variational strategy permits to obtain directly a stable minimal surface, which should produce a corresponding soap film.

\begin{figure}[htbp]
\begin{center}
\includegraphics[width=9cm]{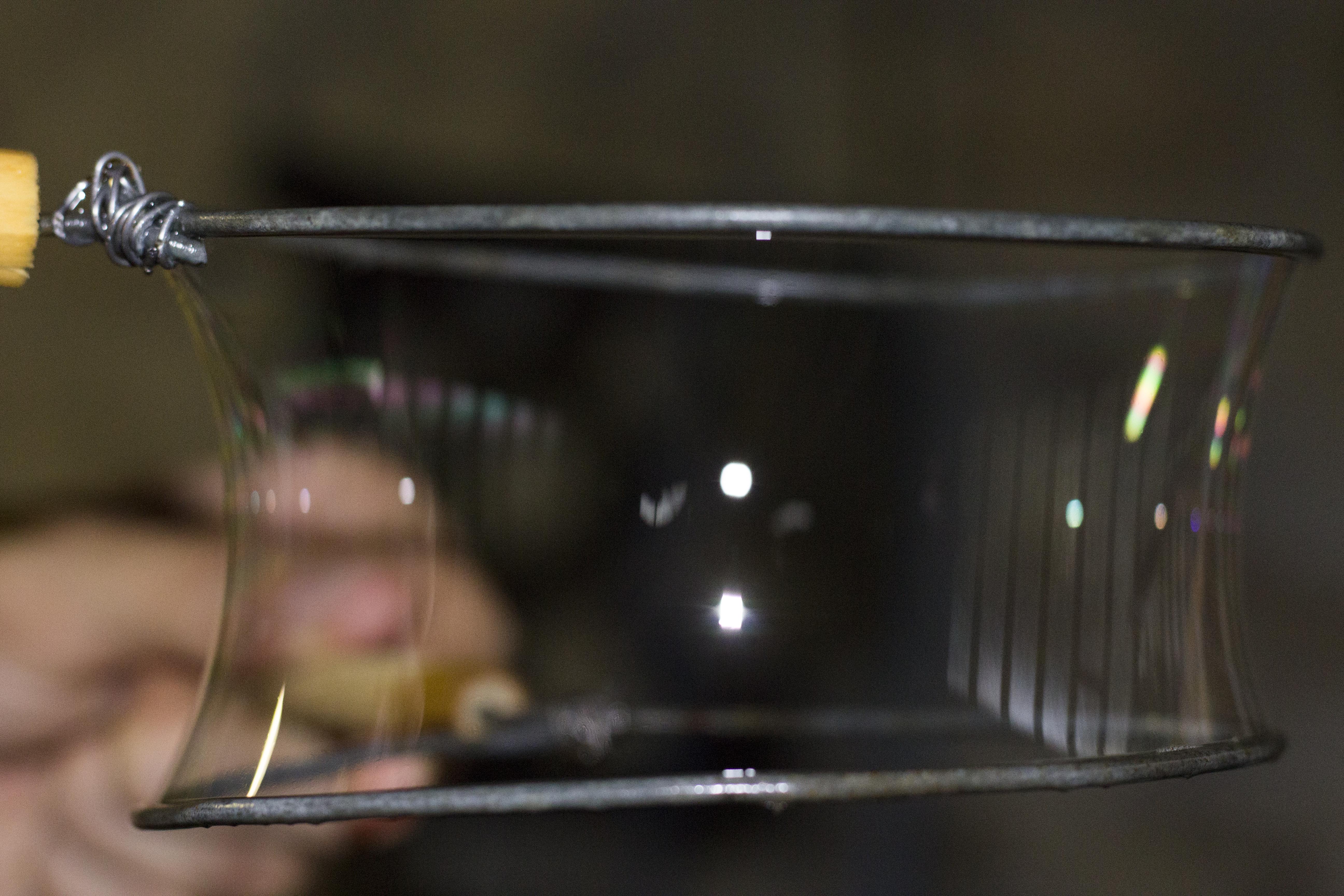} \qquad \includegraphics[width=5.2cm]{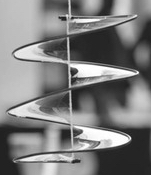}
\caption{The catenoid (on the left) and the helicoid (on the right) as soap films.}\label{sup2}
\end{center}
\end{figure}

However, using some special solutions of the equation $H=0$, one can construct many examples of minimal surfaces. The {\it catenoid}, discovered by Euler in 1744 and obtained as a revolution surface, is the first example of non planar minimal surface. 
The {\it helicoid} is another classical example
as well as the {\it Enneper surface} and the {\it Scherk surface} \cite{do2016differential}. Some of these examples of minimal surfaces can be seen as soap films: the catenoid appears also as a solution of the soap film bounded by two sufficiently close coaxial rings, while the helicoid spans a circular helix (see \cref{sup2}).

\begin{figure}[htbp]
\begin{center}
\includegraphics[width=7cm]{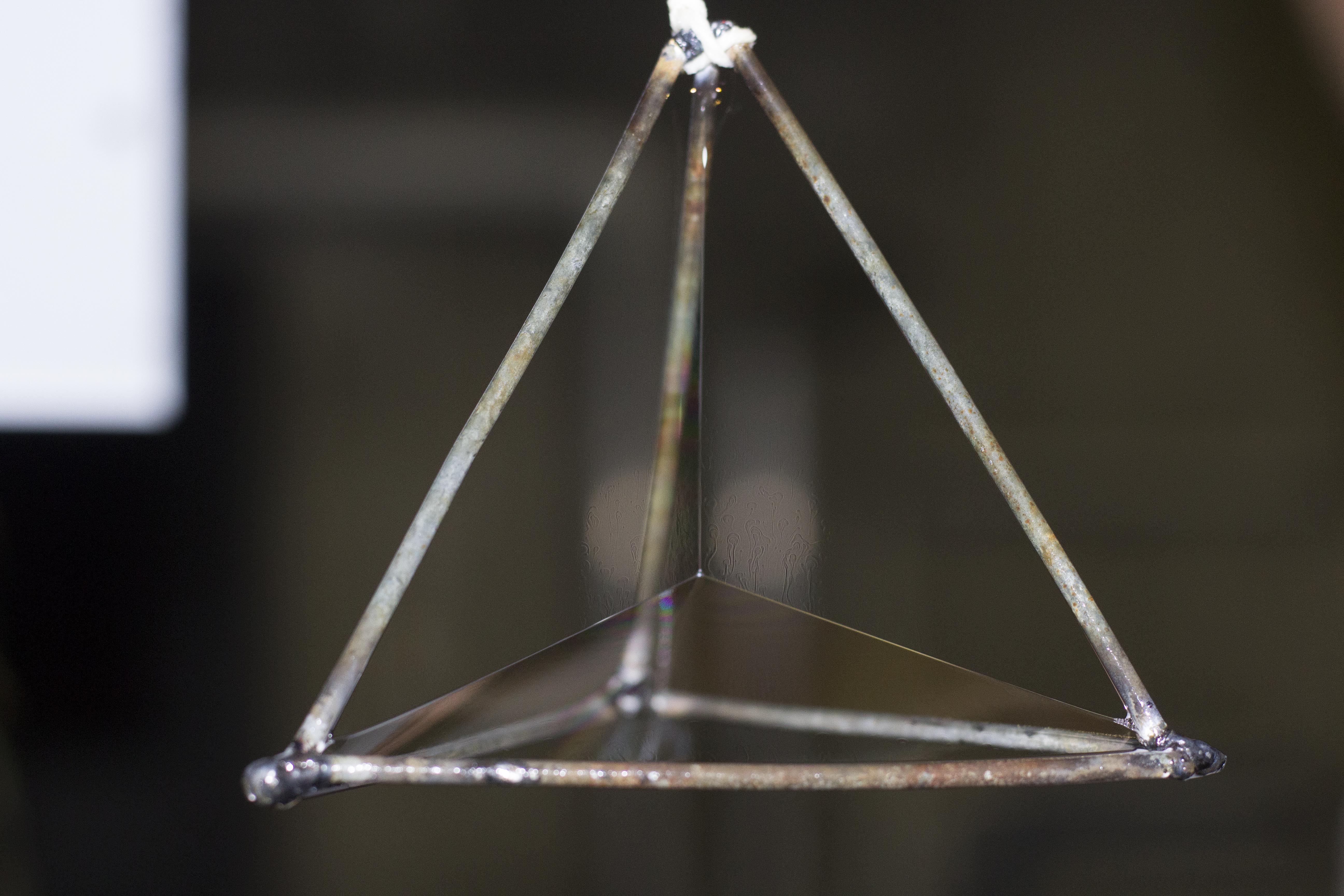} \qquad \includegraphics[width=7cm]{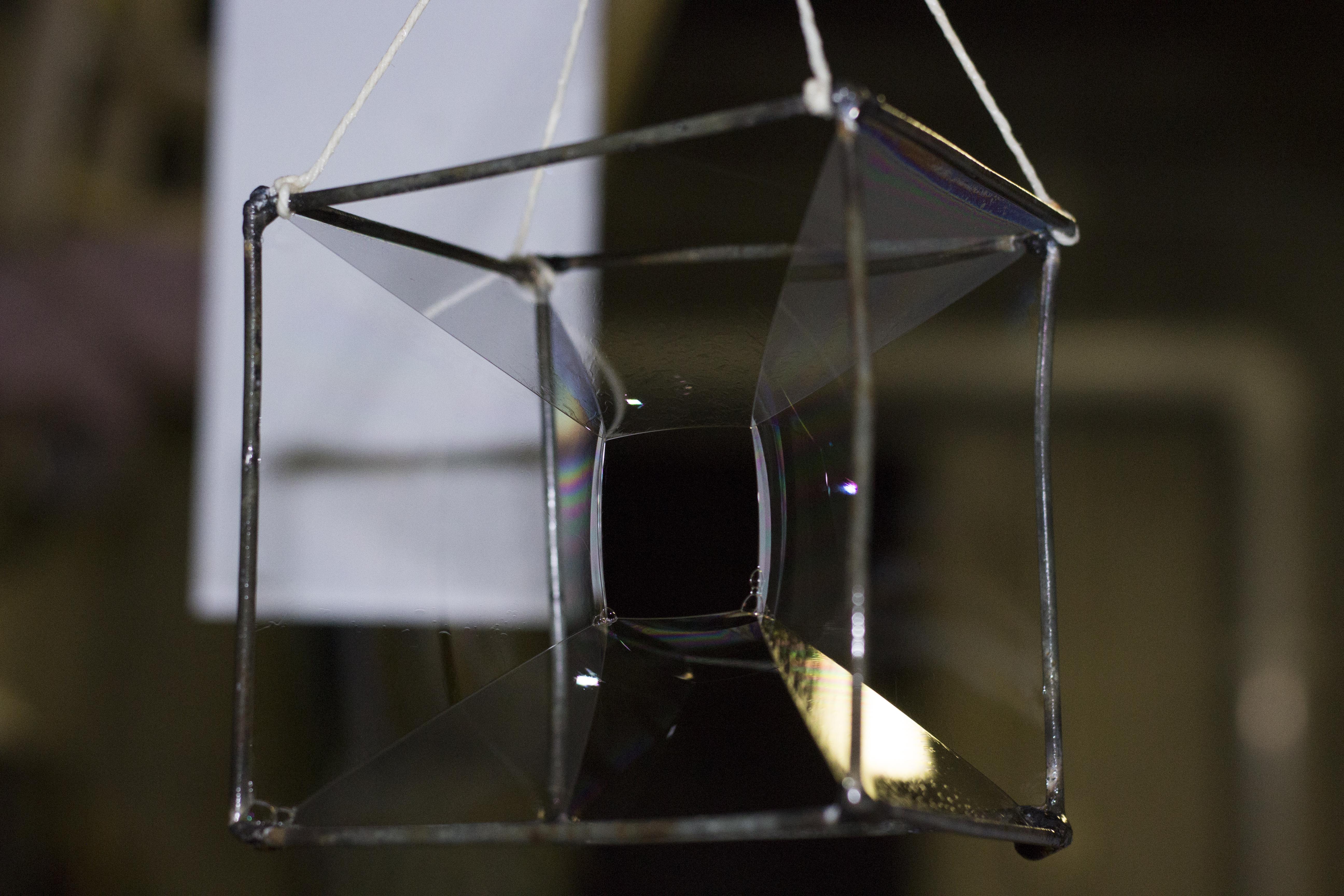}
\caption{The soap film obtained by the edges of a tetrahedron (on the left) and the soap film realized by the edges of a cube (on the right).}\label{other}
\end{center}
\end{figure}

The general formulation of the Plateau problem might be the following one: {\it given a closed curve $\Gamma$ in the space find a surface with minimal area spanning $\Gamma$}. Since it seems that every closed wire spans some soap film, Plateau was convinced that every closed curve with no double points spans a surface which minimizes the area. Moreover, as far as experiments suggest, there are only two kind of singular configurations: the {\it $\mathbb Y$-configuration}, three plane sheets crossing on a line and forming a $120^\circ$ angle, and the {\it $\mathbb T$-configuration}, four lines crossing in a point (called {\it tetrahedrical} point) and forming an angle of approximately $109,47^\circ$. Pictures in \cref{other} show that these singularities may occur.
The $\mathbb Y$ and the $\mathbb T$ singularities are the only conjectured by Plateau. For this reason, the fact that a soap film can only produces $\mathbb Y$ and/or  $\mathbb T$ singularities are known as {\it Plateau laws}. 

In this paper, we will first of all review the main techniques for solving the Plateau problem. From the mathematical point of view, the problem is very difficult and a lot of possible formulations are available. Precisely, in \cref{sec-p} we will briefly mention how the classical solution by Douglas and Rad\'o works, then we will pass to review more recent formulations of the problem in the context of Geometric Measure Theory: sets of finite perimeter, currents, and minimal sets. An important generalization of the Plateau problem is presented in \cref{sec-kp} where the so-called {\it Kirchhoff-Plateau problem} is introduced: the boundary is elastic and can sustain bending and twisting. We will present some recent results for the Kirchhoff-Plateau problem and some of its generalizations. Finally, we conclude in \cref{sec-working_progress} with other results and some future directions of investigations.

\section{Plateau problem: classical and non-classical tools}\label{sec-p}

In order to rigorously state the Plateau problem, we need to clarify three things: what {\it surface} means, what {\it area of a surface} means, and the concept of {\it spanning} a prescribed boundary curve. 

\subsection{Plateau problem for graphs}
The first attempt is to deal with graphs of functions. In this setting, the most elementary problem is to find a {\it global minimal graph}, namely a function $u\colon \R^2\to \R$ whose graph solves the equation ${\bf H}=0$. Any affine function is obviously a solution. Is it the unique solution? This is the celebrated {\it Bernstein problem}, that can be generalized to any dimension: if the graph of a smooth function $\R^n\to \R$ is a minimal surface in $\R^{n+1}$, does this imply that it is an affine function? Bernstein formulated the problem in 1914 and solved it, in the same year, only for $n=2$. 
After that, many mathematicians tried to attack the problem in higher dimensions: we mention Simons who answered positively in $n=6$ and gave an example of locally stable cones in $\R^8$ but without proving that these cones are minimal surfaces on the whole space $\R^8$. Finally, Bombieri, De Giorgi and Giusti showed that Simons cones are indeed minimal surfaces in $\R^n$ for $n\ge 8$. An example is the cone $\{(x,y)\in \R^4 \times \R^4: |x|=|y|\}$.

If $u:\Omega \subset \R^n \to \R$ is a smooth graph, then the Plateau problem reads
\begin{equation}\label{plateaugraphs}
\left\{\begin{array}{ll}
\displaystyle -{\rm div}\frac{\nabla u}{\sqrt{1+|\nabla u|^2}}=0 & \text{on $\Omega$,}\\
u=u_0 & \text{on $\partial\Omega$,}
\end{array}
\right.
\end{equation}
where $\Omega$ is open and bounded in $\R^n$. 
Concerning the solution, we mention Jenkin and Serrin \cite{JS}: if $\partial\Omega$ is of class $C^{2,\alpha}$ for some $\alpha \in (0,1)$, and $u_0\in C^{2,\alpha}(\overline \Omega)$ then \eqref{plateaugraphs} has a solution if and only if the mean curvature of $\partial \Omega$ is everywhere non-negative. On the other hand, if $\partial\Omega$ is of class $C^{2,\alpha}$ for some $\alpha \in (0,1)$, then there exists $\e>0$ such that for every $u_0\in C^{2,\alpha}(\overline \Omega)$ with $\|u_0\|_{2,\alpha}\le \e$ the problem \eqref{plateaugraphs} has a unique solution $u\in C^{2,\alpha}(\overline \Omega)$ \cite{K, Mu}. 
Moreover, \eqref{plateaugraphs} can be stated in a variational way \cite{HS}: it is the Euler-Lagrange equation of the area functional written for graphs like
\begin{equation}\label{farea}
u\mapsto \int_\Omega \sqrt{1+|\nabla u|^2}dx.
\end{equation}
In \cite{HS}, the authors show that if $\Omega$ is convex and $u_0\colon \partial \Omega \to \R$ satisfies the bounded slope condition (that is $u_0$ satisfies a Lipschitz inequality), then the functional \eqref{farea} has a unique minimizer among all Lipschitz functions with $u=u_0$ on $\partial\Omega$.

\subsection{Disc-type solutions}

The first rigorous solution to the Plateau problem is due to Douglas \cite{Douglas} and Rad\'o \cite{Rado}, who independently developed an argument which works only for $2$-dimensional surfaces in $\R^3$ and in codimension $1$. We also mention simplifications in the proof of Courant, Tonelli and Dierkes \cite{DHS}. 
The basic idea is to look at smooth parametrizations ${\bf X}\colon D \to \R^3$ where $D=\{(u,v)\in \R^2 : u^2+v^2<1\}$ is the disc and the trace of ${\bf X}$ on $\partial D$ is a smooth parametrization of a prescribed Jordan curve $\Gamma$ in $\R^3$. Thus, the Plateau problem reads
\[
{\bf A}({\bf X})=\int_D|\partial_u{\bf X}\times\partial_v{\bf X}|\,dudv.
\]
In order to apply the Direct Method of the Calculus of Variations, one immediately notices that ${\bf A}$ is weakly lower semicontinuous since the map $(u,v)\mapsto |\partial_u{\bf X}(u,v)\times\partial_v{\bf X}(u,v)|$ is a convex function of the determinants of the $2\times 2$ minors of $\nabla {\bf X}$ (what is called {\it polyconvex} function \cite{D}). 
Concerning compactness, the set $\{{\bf X}: {\bf A}({\bf X})\le c\}$ is not bounded in any reasonable Sobolev norm since the area functional is invariant under reparametrization. However, considering conformal coordinates, it is easy to see that 
\[
{\bf A}({\bf X})\le \frac{1}{2}\int_D |\nabla {\bf X}|^2\,dudv={\bf D}({\bf X}),
\]
and the equality holds true if and only if ${\bf X}$ is conformal, that is $|{\bf X}_u|=|{\bf X}_v|$ and ${\bf X}_u\cdot {\bf X}_v=0$. 
This suggests to minimize directly the Dirichlet functional ${\bf D}$, which is not invariant under reparametrization, among all ${\bf X}\in W^{1,2}(D;\R^3)$ such that ${\bf X}_{|_{\partial D}}$ is a reparametrization of $\Gamma$. 
Precisely, there exists a minimizer ${\bf X}_0\in C^0(\overline D;\R^3)$ which is harmonic on $D$ and ${\bf X}_0$ is conformal, hence ${\bf D}({\bf X}_0)={\bf A}({\bf X}_0)$. 
It is left to show that every minimizer ${\bf X}_0$ of ${\bf D}$ satisfying ${\bf D}({\bf X}_0)={\bf A}({\bf X}_0)$ is a minimizer for the area functional. Obvioulsly, since for any admissible ${\bf X}$ it holds ${\bf A}({\bf X})\le {\bf D}({\bf X})$, then
\[
\inf_{{\bf X}}{\bf A}\le \inf_{{\bf X}}{\bf D}.
\]
To have the equality one can apply the $\e$-conformal mappings Lemma due to Morrey \cite{Mo}: if ${\bf X}\in C^0(\overline D;\R^3) \cap W^{1,2}(D;\R^3)$ then for any $\e > 0$ there exists a homeomorphism $\tau_\e \colon \overline D\to \overline D$ of class $W^{1,2}$ such that ${\bf D}({\bf X}\circ \tau_\e) \le {\bf A}({\bf X}) + \e$. 
Finally, ${\bf X}_0$ produces a regular surface, namely $\partial_u{\bf X}\times \partial_v{\bf X}\ne 0$ everywhere, indeed if $\Gamma$ is an analytical Jordan curve and if its total curvature does not exceed $4\pi$ then any disc-type solution of Plateau problem is a regular minimal surface \cite{N}.
\begin{figure}[htbp]
\begin{center}
\includegraphics[width=7.5cm]{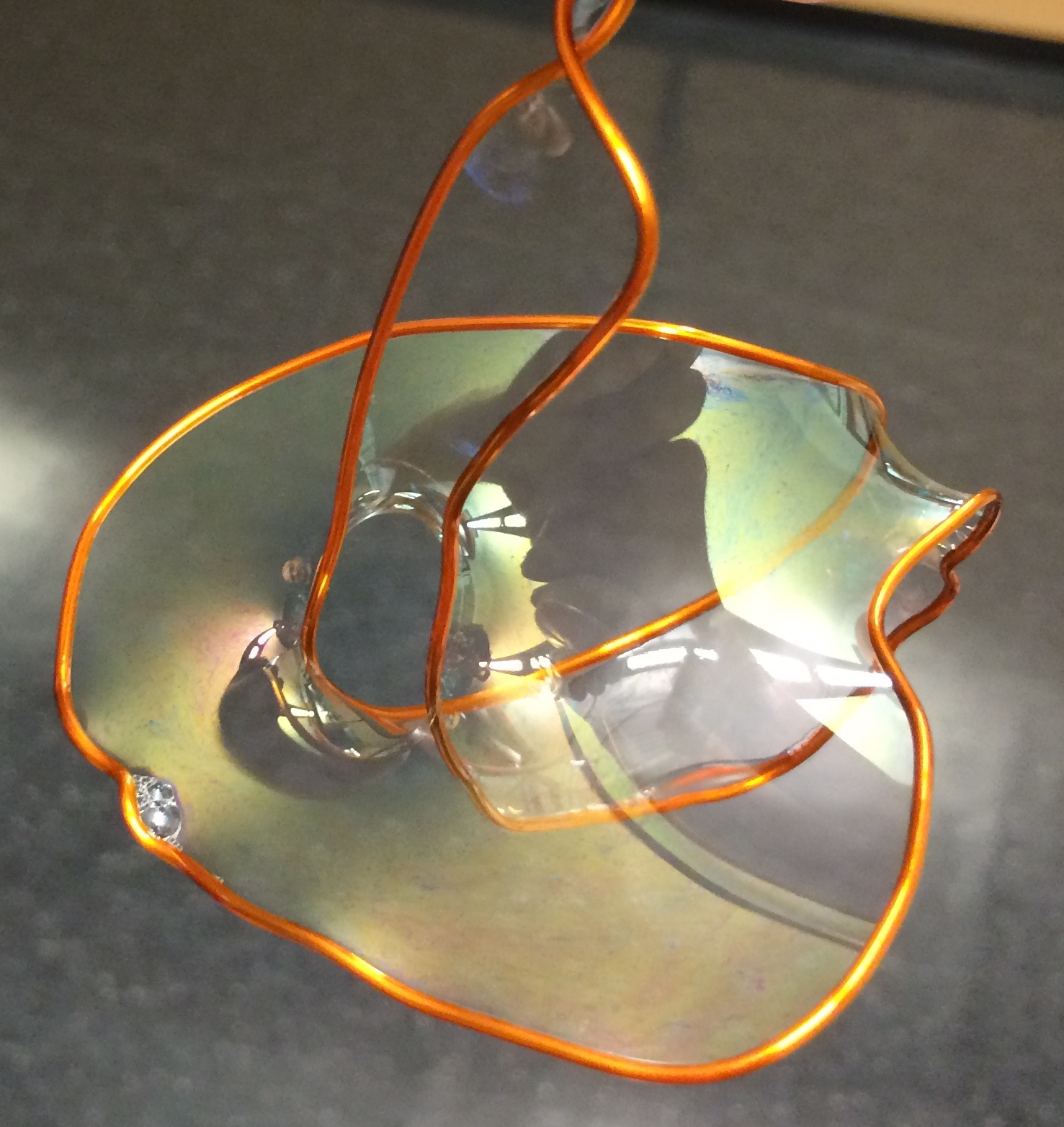}
\caption{The area minimizing soap film spanning a disc wants to be embedded.}\label{fig3}
\end{center}
\end{figure}
Moreover, disc-type minimal surfaces cannot produce singularities in the interior \cite{G}: they are immersed surfaces and not embedded, without self-intersections (see for instance the embedded soap film solution in \cref{fig3}), not providing a good model for soap films.

\subsection{Distributional approaches}

Concerning distributional approaches to solve the Plateau problem, we describe two approaches: sets of finite perimeter and currents.

We refer to the book by Ambrosio-Fusco-Pallara \cite{AFP} or to the monograph by Maggi \cite{Maggi} for details on the theory of finite perimeter sets. Let $E\subset \R^n$ be a Borel set with finite Lebesgue measure and let $\Omega\subseteq \R^n$ be open. We say that $E$ is a {\it finite perimeter set in $\Omega$} if 
\[
\mathscr P(E;\Omega)=\sup\left\{\int_{E\cap \Omega}  {\rm div}\,\phi\,dx : \phi\in C^\infty_c(\R^n;\R^n),\,\,\|\phi\|_\infty\le 1\right\}<+\infty.
\]
The quantity $\mathscr P(E;\Omega)$ is the {\it perimeter of $E$ in $\Omega$}. If $E$ is a bounded and open with smooth boundary, then $E$ has finite perimeter in $\Omega$ and $\mathscr P(E;\Omega)=\mathscr H^{n-1}(\partial E\cap \Omega)$. By duality, we immediately can say that the map $E\mapsto \mathscr P(E;\Omega)$ is lower semicontinuous with respect to the $L^1$-convergence of characteristic function f sets. Moreover, if $\mathscr P(E_h)$ is bounded and $E_h$ are contained in a given ball, then, up to a subsequence, $E_h \stackrel{L^1}{\to}E$ and $E$ has finite perimeter.

Thus, the Plateau problem in this setting reads as follows:
\begin{equation}
    \label{eq:minimi_perimetro}
    \inf\left\{\mathscr P(E) : \text{$E$ has finite perimeter in $\Omega$ satisfying } \mathscr{L}^n((E\setminus\Omega)\mathbin{\scriptstyle\bigtriangleup} E_0)=0\right\},
\end{equation}
where $E_0 \subset \R^n \setminus \Omega$ be such that $\partial E_0 \cap \partial \Omega = \Sigma_0 \subset \partial \Omega$ assigned. In particular, the Direct Method of the Calculus of Variations can be successfully applied and the problem \eqref{eq:minimi_perimetro} has a minimal solution.

Another distributional approach to the Plateau problem is the use of the theory of currents. The notion of current dates back to De Rham, while the variational and geometrical approach used today is mainly due to Federer and Fleming \cite{F, S}.

Let $\mathscr D^d(\R^n)$ be the set of $d$-forms on $\R^n$ with compact support. The space of {\it $d$-currents on $\R^n$}, denoted by $\mathscr D_d(\R^n)$, is the topological dual space of $\mathscr D^d(\R^n)$. 
Then, any $d$-dimensional smooth oriented surface $S$ in $\R^n$ is an example of a current: $T_S\in \mathscr D_d(\R^n)$ is defined as follows
\[
\langle T_S,\omega\rangle=\int_S\omega, \qquad \forall \omega \in \mathscr D^d(\R^n).
\]
Also the boundary of a current can be defined via the Stokes' formula: if $T\in \mathscr D_d(\R^n)$, then $\partial T\in  \mathscr D_{d-1}(\R^n)$ is the {\it boundary of $T$} and it is given by 
\[
\langle \partial T,\omega\rangle=\langle T,d\omega \rangle. 
\]
Moreover, by Stokes' formula and for smooth oriented surfaces, $\partial T_S=T_{\partial S}$.

To state the Plateau problem, the concept of the {\it mass} of a current $T\in\mathscr D_d(\R^n)$ is defined as
\[
\mathbb M(T)=\sup_{||\omega(x)||\le 1}\langle T,\omega\rangle,
\]
where $\|\omega(x)\|$ is a suitable notion of norm for $d$-forms. It turns out that $\mathbb M$ is lower semicontinuous with respect to the weak convergence of currents. For smooth oriented surfaces, it holds $\mathbb M(T_S)=\mathscr H^d(S)$. 

Since the space of currents is too large, a subspace must be introduced: $T\in\mathscr D_d(\R^n)$ is a {\it $d$-rectifiable current with integer multiplicity} if there exist:
\begin{itemize}
\item[\rm(a)] a $d$-rectifiable set $E$ in $\R^n$,
\item[\rm(b)] an {\it orientation $\tau$ on $E$}, namely a Borel map that to $\mathscr H^d$-a.e.\,$x\in E$ assigns a unit simple $d$-vector $\tau(x)$ which spans $T_x E$,
\item[\rm(c)] a {\it multiplicity function}, that is a $\mathscr H^d$-summable function $m\colon E\to \N$,
\end{itemize}
such that $T$ can be defined as follows
\[
\langle T,\omega\rangle=\int_E\langle\omega(x),\tau(x)\rangle\, m(x)\,d\mathscr H^d(x), \qquad  \forall \omega \in \mathscr D^d(\R^n).
\]
A current of this type is denoted by $[E,\tau,m]$. If $S$ is a smooth $d$-dimensional surface oriented by $\tau$ then $T_S=[S,\tau,1]$. 

Finally, a current $T\in\mathscr D_d(\R^n)$ is said to be a {\it $d$-integral current} if both $T$ and $\partial T$ are rectifiable currents with integer multiplicity. In this final class, a compactness theorem holds true (Federer-Fleming Compactness Theorem): if $(T_h)$ is a sequence of integral $d$-currents with $\mathbb M(T_h)+\mathbb M(\partial T_h)$ bounded then, up to a subsequence, $T_h\to T$ in $\mathscr D_d(\R^n)$, $\partial T_h\to \partial T$ where $T$ is an integral current. Moreover,  
\[
\mathbb M(T)\le \liminf_h \mathbb M(T_h), \qquad \mathbb M(\partial T)\le \liminf_h \mathbb M(\partial T_h).
\]
The Plateau problem in terms of integral currents can be stated as follows: let $T_0$ be a given integral $d$-current on $\R^n$; find a minimizer of $\mathbb M(T)$ among all currents $d$ integral currents $T$ with $\partial T=\partial T_0$.

\begin{figure}[htbp]
\begin{center}
\includegraphics[width=9.5cm]{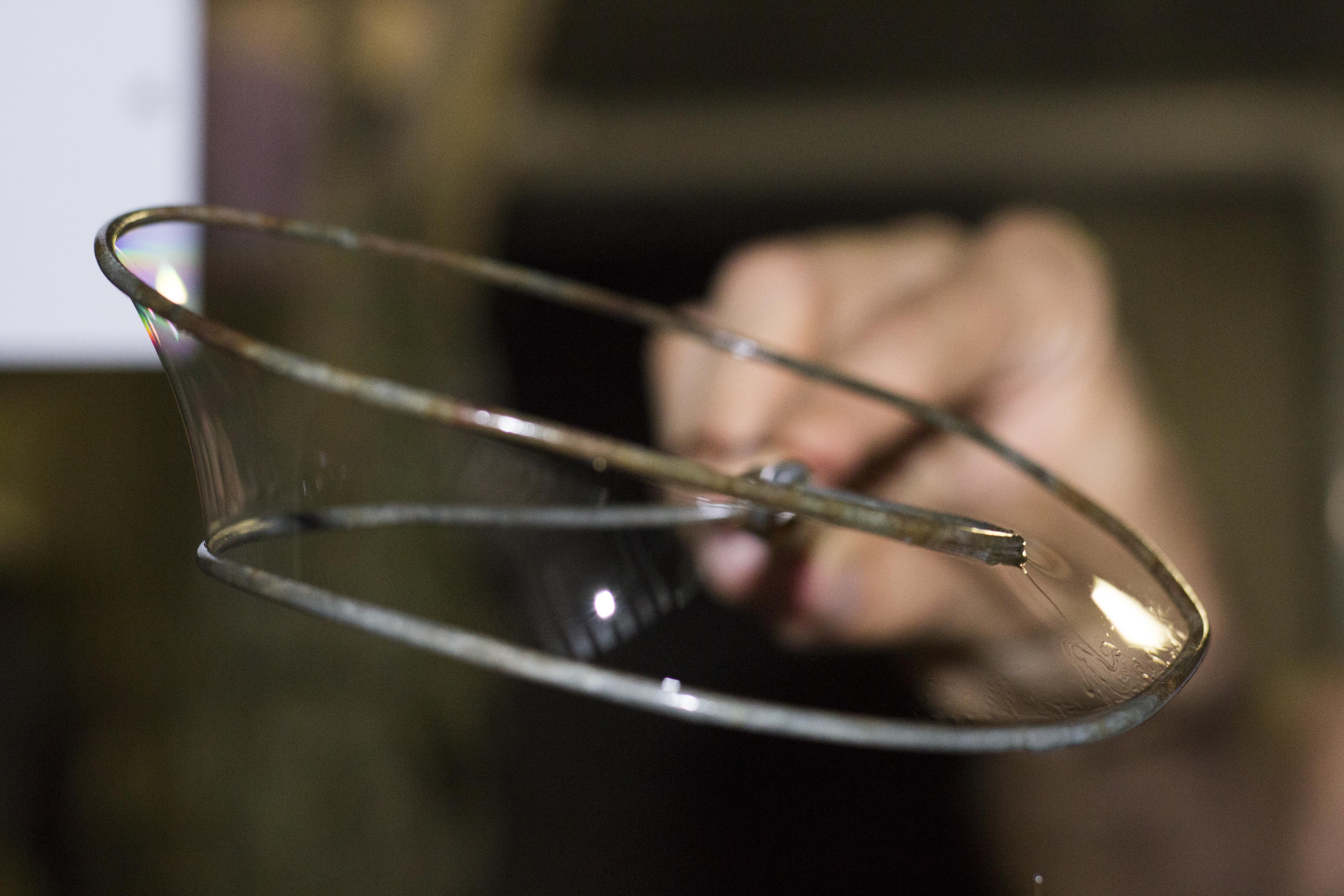}
\caption{A M\"obius strip-like soap film}\label{mobius}
\end{center}
\end{figure}

Notice that this approach has some limitations:
\begin{enumerate}
    \item the current solution to the Plateau problem can have multiplicity different from $1$. Indeed, the right object to minimize should be the {\it size} of a current defined as
\[
\mathbb S\left([E,\tau,m]\right)=\mathscr H^d\left(\{x\in E : m(x)\ne 0\}\right).
\]
However, for $\mathbb S$ a compactness theorem does not hold true.
\item any discontinuity on the orientation produces an nonphysical boundary. A possibility to overcome this issue is to produce non-orientable soap films (see \cref{mobius}) and mathematically to deal with rectifiable currents {\it modulo $\nu$}, where $\nu\ge 2$ is an integer or using the theory of varifold, for definitions and details see \cite{A, H}.
\end{enumerate}

Finally, we would like to say that in both cases a minimizer produces actually a soap film. Thus, a regularity theory has been developed. 
For the set of finite perimeter the regularity is obtained for a set $E$ in $\R^n$ which minimizes the perimeter with respect to all possible compactly supported perturbations. In this case it is possible to prove that $\partial E\setminus S$ is smooth, where $S$ is the closed set of singularities and
\begin{itemize}
\item[\rm(a)] if $2\le n\le 7$ then $S$ is empty and $\partial E$ is analytical;
\item[\rm(b)] if $n = 8$ then $S$ has no accumulation points in $E$;
\item[\rm(c)] if $n\ge 9$ then $\mathscr H^d(S)=0$ for every $d>n-8$.  
\end{itemize}
A similar regularity result holds for mass-minimizing currents: if $T$ is a mass-minimizing 1-integral current in $\R^2$ then the ``interior part" of $T$ (the part of $T$ which is not in the boundary of $T$) is made of disjoint line segments. Moreover, similar to the set of finite perimeter case, if $2\le n\le 7$ then the interior part of any mass-minimizing $(n-1)$-integral current in $\R^n$ is a smooth embedded hypersurface: in \cref{fig3} the soap film solution corresponds to an embedded solution, which actually should be the mass-minimizing integral current.
When $n>7$, the Simons cone $C=\{(x,y)\in \R^4\times \R^4 : |x|=|y|\}$ is an area-minimizer current developing a singularity in the origin \cite{bdgg}.
Unfortunately, since in lower dimension, especially the physical one $n =3$, both the set of finite perimeter approach and the current one do not develop singularities, they do not provide a good model to study soap films.


\subsection{Almgren minimal sets approach and Taylor regularity}\label{sec:minimalsets}

Minimal sets, introduced by Almgren in \cite{A1}, represent the best model for soap films. 

Let us recall the definition of minimal set. Let $S\subset \R^n$ be a closed set and $A \subset \R^n$ be an open set. We say that $S$ is a {\it $d$-dimensional minimal set in $A$}, briefly {\it minimal set}, if for any closed ball $B\subset A$ and for every Lipschitz function $\varphi\colon \R^n\to \R^n$ with $\varphi_{|_{\R^n\setminus C}}=id$ and with $\varphi(C)\subset C$ we have $\mathscr H^d(S)\le \mathscr H^d(\varphi(S))$. In 1976, Taylor \cite{T} proved that $2$-dimensional minimal sets in $\R^3$ may have singularities and these are exactly the ones produced by soap films and observed by Plateau in his experiments. More precisely, there are two kind of singularities (see \cref{other}):
\begin{itemize}
\item[\rm(a)] the so-called {\it $\mathbb Y$-configuration}: three sheets crossing on a line and forming a $120^\circ$ angle;
\item[\rm(b)] the so-called {\it $\mathbb T$-configuration}: four lines crossing in a point forming a $109,47^\circ$ angle.
\end{itemize}

\begin{figure}[htbp]
\begin{center}
\includegraphics[width=16.5cm]{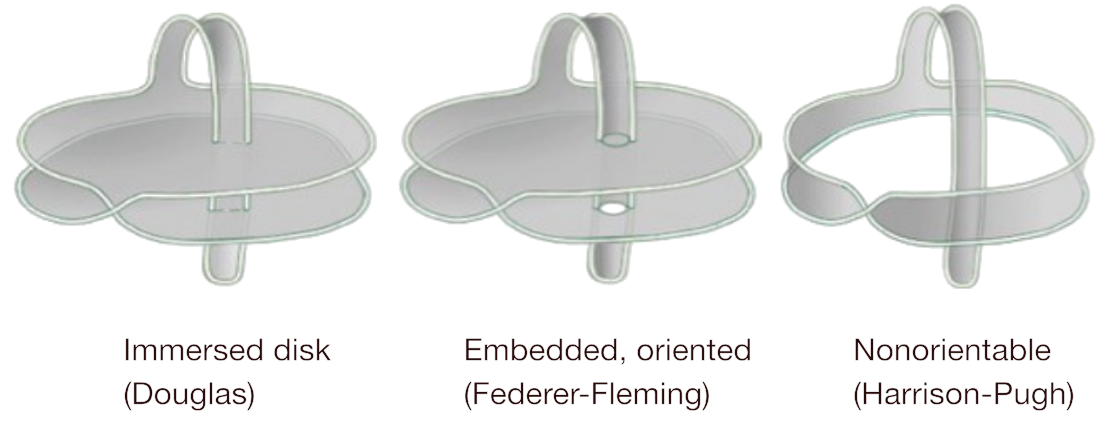} 
\caption{Three different solutions for the same wire (courtesy of J.\,Harrison \cite{H}).}
\label{h1}
\end{center}
\end{figure}

In order to state a Plateau problem in this framework the main  difficulty stems from the notion of boundary. Recently, a suitable theory has been developed and some existence results have been proved. The main idea has been introduced by Harrison \cite{Ha, HP}. The approach by Harrison and Pugh is based on differential chains and it permits to represent all types of observed soap films as well as immersed surfaces of various genus types, both orientable and nonorientable, see \cref{h1}.

Later on, De Lellis, Ghiraldin, and Maggi reformulated the concept of spanning in a more Geometric Measure Theory setting  \cite{DGM}: let $n\ge 3$ and let $H$ be a closed subset of $\R^n$. Let 
\[
\mathscr C_H=\{\gamma \colon \mathbb S^1 \to \R^n\setminus H\,\,\textrm{smooth embedding of $\mathbb S^1$ into $\R^n$}\}. 
\]
Fix $\mathscr C\subset \mathscr C_H$ be a closed subset by homotopy and let $K$ be a relatively closed set in $\R^n\setminus H$. The set $K$ is said to be a {\it $\mathscr C$-spanning set of $H$} if 
\begin{equation}
\label{e:spanning}
    K \cap \gamma(\mathbb S^1) \ne \emptyset,\qquad  \forall \gamma \in \mathscr C.
\end{equation}
Let us denote by $\mathscr F(H,\mathscr C)$ the class of all relatively closed sets in $\R^n\setminus H$ which are $\mathscr C$-spanning sets of $H$. If there exists $K\in \mathscr F(H,\mathscr C)$ such that $\mathscr H^{n-1}(K)<+\infty$, then the problem 
\[
\min_{K\in \mathscr F(H,\mathscr C)}\mathscr H^{n-1}(K)
\]
has a solution which is a $(n-1)$-dimensional minimal set in $\R^n\setminus H$. 

This approach furnishes a good answer to the Plateau problem: when $H$ is a Jordan curve in $\R^3$ the spanning condition corresponds to the fact that the soap film $K$ wets entirely the curve $H$. There exists a minimal set $K$ in $\R^3\setminus H$ that spans $H$. Therefore, ``the boundary of $K$ is $H$'' and, by Taylor's result, $K$ can develop Plateau's type singularities. 

\section{The Kirchhoff-Plateau problem}\label{sec-kp}

A recent generalization of the Plateau problem consists in the situation in which the boundary is not rigid but it is given by a flexible manifold. The first results were given by Bernatzky \cite{B} and Bernatzky and Ye \cite{BY} who proved existence in the framework of currents. The first formulation of the Kirchhoff-Plateau problem was given by Giusteri, Franceschini and Fried \cite{GFF}, where the boundary of the soap film lies on a $3$-dimensional elastic rod; more precisely, stability of equilibrium configurations is analyzed. We also mention \cite{CheFri14, BirFri14, BirFri15, HoaFri16}.
The first rigorous existence result for the Kirchhoff-Plateau problem has been provided by  Fried, Giusteri and Lussardi \cite{GLF} where the energy functional to be minimized is composed by the elastic energy of the rod, the weight of the rod and the area of the soap film spanned by the rod. Further details and relative bibliography can be found therein; we refer also to \cite{BLM2, BLM1, BLM3, DL, bevilacqua2024effects}.

\subsection{The bounding loop}

To model the boundary manifold, in \cite{GLF} the theory of Kirchhoff rods has been implemented (see for instance the book of Antman \cite{Antman2005}). A $3D$-rod is completely described by its midline curve and a family of two-dimensional {\it material cross-section} attached to each point of the midline. Moreover, in order to encode  how the cross-sections are ``appended'' to the midline, a family of {\it material frames} completes the framework. Here, it is also assumed that the material cross-section lies in the plane orthogonal to the midline at any point of the midline, namely that the rod is {\it unshearable}, and that its midline is {\it inextensible}, denoting wit $L>0$ its length.
Under these assumptions, the final shape of the rod is uniquely determined by assigning the cross-sections and three scalar fields: the {\it flexural densities} $\kappa_1$ and $\kappa_2$ and a {\it twist density} $\omega$ as we illustrate. Fix the clamping point $\vect x_0 \in \R^3$ and fix $\vect t_0,\vect d_0 \in \R^3$ unit orthogonal vectors. Let $p>1$ and $V=L^p(0,L)\times\R^3\times\R^3\times\R^3$, then $\vect w \in V$ and it is given by
\[
\vect w=((\kappa_1,\kappa_2,\omega),\vect x_0,\vect t_0,\vect d_0).
\]
Starting from $\vect w$, we can reconstruct the midline $\vect x$ and a director field $\vect d$ as the unique solutions of the following system of ordinary differential equations (for a graphical representation see \cref{rod}) 
\begin{equation}\label{eq:Cauchy}
\left\{\begin{aligned}
&{\vect x}'=\vect t,\\
&{\vect t}'=\kappa_1\vect d+\kappa_2\vect t\times \vect d,\\
&{\vect d}'=\omega\vect t\times \vect d-\kappa_1\vect t,
\end{aligned}\right.
\end{equation}
supplemented by the initial conditions
\begin{equation}\label{eq:CauchyIC}
\left\{\begin{aligned}
&\vect x(0)=\vect x_0,\\
&\vect t(0)=\vect t_0,\\
&\vect d(0)=\vect d_0.
\end{aligned}\right.
\end{equation}
\begin{figure}[htbp]
\begin{center}
\includegraphics[width=0.95\textwidth]{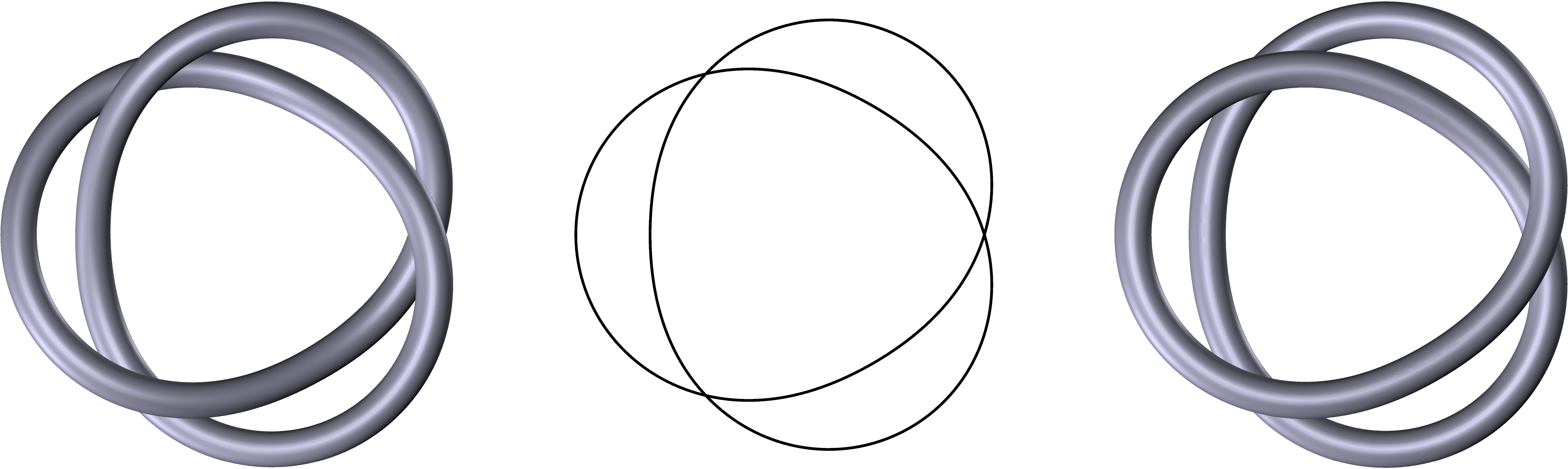}
\end{center}
\caption{The trefoil knot on the left and the unknot structure on the right  are distinct objects when the cross-sections have a non-vanishing thickness, even in the presence of self-contact. The distinction between a knot and an unknot is lost in the vanishing-thickness limit.}
\label{fig:thickness}
\end{figure}
A classical result of Carath\'eodory \cite{Hartman1982} ensures that \eqref{eq:Cauchy}--\eqref{eq:CauchyIC} has a unique solution:
$$
\vect x \in W^{2,p}((0,L);\R^3) \qquad \hbox{ and } \qquad \vect d \in W^{1,p}((0,L);\R^3).
$$
Moreover, since $\vect t_0,\vect d_0$ are unit orthogonal vectors, we can reconstruct the moving material frame given by $\{(\vect t(s),\vect d(s),\vect t(s)\times \vect d(s)):s\in[0,L]\}$. In particular, the midline $\vect x$ is parametrized by the arc-length. In addition, the problem is equipped with some constraints.
First of all, in order to have a closed and smooth midline we require that 
\begin{equation}\label{eq:closure}
\vect x(L)=\vect x(0), \qquad \hbox{ and } \qquad \vect t(L)=\vect t(0).
\end{equation}
We want also to encode the knot type of midline. To do that we simply fix a continuous map $\ell\colon [0,L]\to\R^3$ with $\ell(L)=\ell(0)$ and we ask that 
\begin{equation}\label{eq:knot}
\vect x\simeq \ell,
\end{equation}
where $\simeq$ is the isotopy equivalence relation in the sense of the theory of knots. We point out that a non-vanishing cross-sectional thickness is crucial for distinguishing knot types in the presence of self-contact, see \cref{fig:thickness}. 
Next, we discuss how to reconstruct the shape of the rod. The material cross-section at each $s \in [0,L]$ is given by a compact and simply connected set $A(s)\subset \R^2$ which contains the origin.
The corresponding  rod can then be described as the set $\vect p[\vect w](\Omega)$, where 
\begin{equation*}
\Omega=\big\{(s,\zeta_1,\zeta_2):s\in[0,L]\text{ and }(\zeta_1,\zeta_2)\in A(s)\big\}
\end{equation*}
and $\vect p[\vect w]$ is given by 
\[
\vect p[\vect w](s,\zeta_1,\zeta_2)=\vect x(s)+\zeta_1\vect d(s)+\zeta_2\vect t(s)\times \vect d(s).
\]
In what follows, $\Lambda[\vect w]$ will stand for the set $\vect p[\vect w](\Omega)$. 
We define the elastic energy of the rod as 
\[
E_\mathrm{sh}(\vect w)=\int_0^Lf(w_1(s),s)\,ds
\]
where $f\colon \R^3\times[0,L]\to\R\cup\{+\infty\}$ satisfies the following conditions:
\begin{itemize}
\item[\rm(a)] $f(\cdot,s)$ is continuous and convex for any $s$ in $[0,L]$;
\item[\rm(b)] $f(a,s)$ is bounded from below;
\item[\rm(c)] $f(a,\cdot)$ is measurable for any $a \in \R^3$.
\item[\rm(d)] $f(a,s)\geq c_1|a|^p+c_2$ for some $c_1,c_2\in \R$ with $c_1>0$.
\end{itemize}
\begin{figure}[htbp]
\begin{center}
\includegraphics[width=11cm]{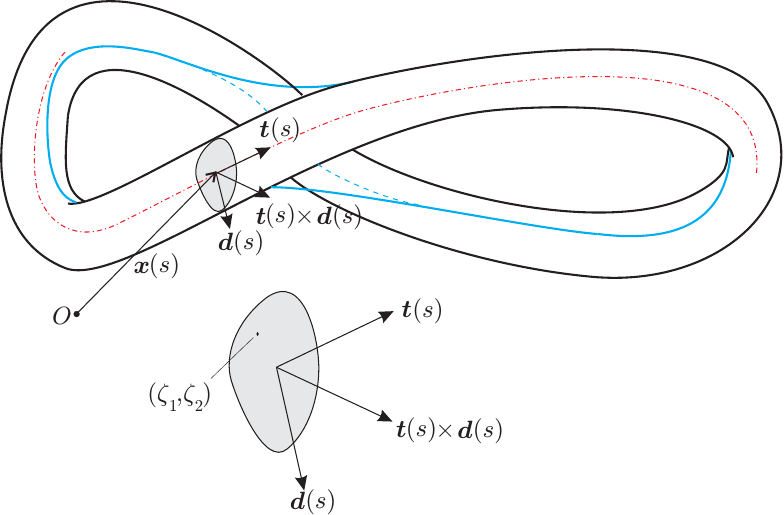} 
\caption{The shape of the rod constructed by a moving frame.}\label{rod}
\end{center}
\end{figure}
Under these assumptions the functional $E_\mathrm{sh}$ is coercive and lower semicontinuous with respect to the weak topology of $V$. We now discuss other necessary physical constraints. 
First of all, we have to specify how we glue the last cross-section to the first one when we close the rod. More precisely, we have to prescribe how many times the ends of the rod are twisted before being glued together. 
Thus, for a small parameter $\e >0 $, the curve $\vect x+\varepsilon \vect d$ remains inside the rod $\Lambda[\vect w]$. Up to add a straigth line-segment we can assume that the curve $\vect x+\varepsilon \vect d$ is closed. We ask that the 
linking number between the midline $\vect x$ and the curve $\vect x+\varepsilon \vect d$ is a prescribed integer number $z$, i.e.
\begin{equation}\label{link_constraint}
{\rm Link}(\vect x,\vect x+\varepsilon \vect d)=z.
\end{equation}
To complete the global gluing condition, see \cref{fig:glue}, we also fix the angle between $\vect d_0$ and $\vect d(L)$ as
\begin{equation}\label{d_constraint}
\textrm{angle\,$(\vect d_0,\vect d(L))$ is fixed}.
\end{equation}
\begin{figure}[htbp]
		\centering
		\includegraphics[width=10cm]{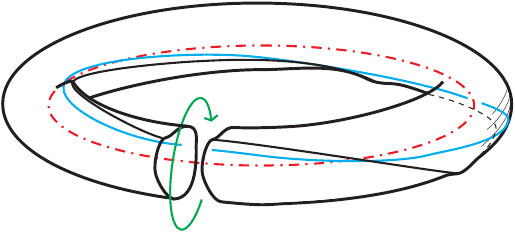}
	\caption{The gluing of the rod: the curve in blue is close to the midline in red.}
 \label{fig:glue}
\end{figure}

\begin{figure}[htbp]
\begin{center}
\includegraphics[width=0.7\textwidth]{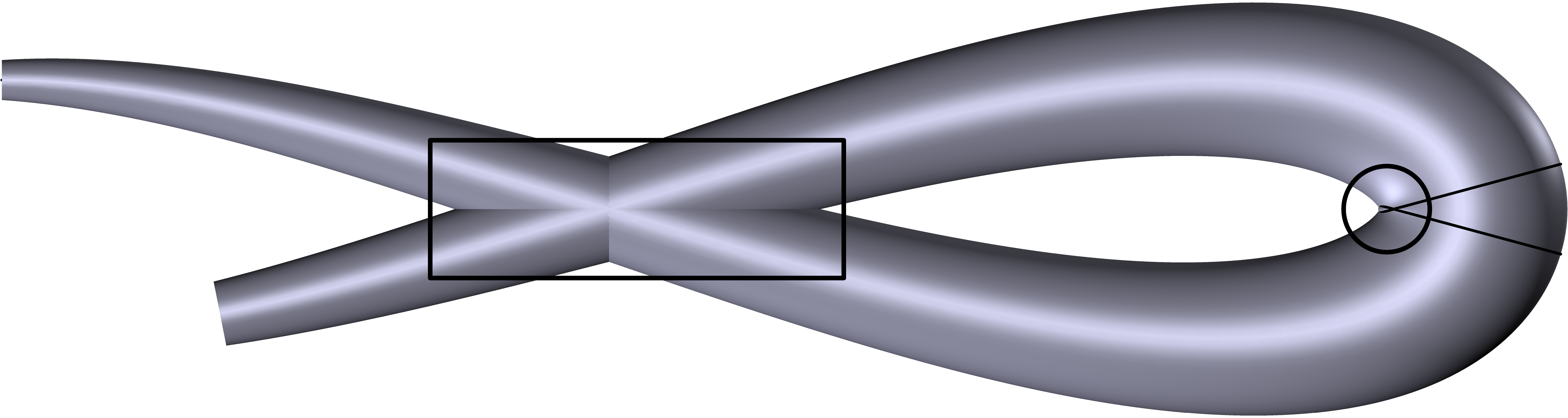}\hspace{.2cm}
\includegraphics[width=0.18\textwidth]{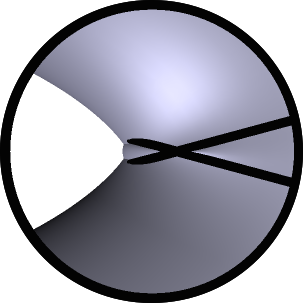}
\end{center}
\caption{Lost of global and local injectivity.}
\label{fig:injectivity}
\end{figure}

\noindent We finally have to discuss the non-interpenetration of matter, see \cref{fig:injectivity}. In order to guarantee that the map $\vect p$ is globally injective on the interior part of $\Omega$ we have to assume two conditions. The first one is the {\it local non-interpenetration constraint}, which we employ adding to the energy of the loop the term
\[
E_{\mathrm{ni}}(\vect w)=
\begin{cases}
0 & \text{if }\vect w\in N\,,\\
+\infty & \text{if }\vect w \in V\setminus N\
\end{cases}
\]
where 
\[
N=\left\{\vect w\in V : \max_{(\zeta_1,\zeta_2)\in A(s)}\big(\zeta_1\kappa_2(s)-\zeta_2\kappa_1(s)\big) \leq 1,\,\text{a.e.\,$s\in (0,L)$}\right\}.
\]
The penalization of $E_{\mathrm{ni}}$ can be seen as a sort of relaxation of the orientation-preservation of $\vect p[\vect w]$. Besides, the global injectivity follows from the Ciarlet-Ne\v cas condition 
\begin{equation}\label{eq:glob-inj}
\int_\Omega\det D\vect p[\vect w](s,\zeta_1,\zeta_2)\,dsd\zeta_1d\zeta_2\leq \mathscr L^3(\Lambda[\vect w]).
\end{equation} 
It remains to add the effects of the weight of the rod: we consider the potential energy
\[
E_{g}(\vect w)=-\int_\Omega \rho(s,\zeta_1,\zeta_2)\,\vect g\cdot \vect p(s,\zeta_1,\zeta_2)
 \, dsd\zeta_1d\zeta_2,
\]
where $\rho>0$ is the mass density and $\vect g$ is the constant acceleration of gravity. 

The final form of the loop energy reads as 
\[
E_\mathrm{loop}=E_\mathrm{sh}+E_{\mathrm{ni}}+E_{g}.
\]

\subsection{The Kirchhoff--Plateau problem}
As in the classical Plateau problem, we model the liquid film by a two-dimensional object $K$, but we want to keep track of the fact that it is reminiscent of two adhering surfactant leaflets. Then, we define the energy of the liquid film as
\[
E_\mathrm{film}(K)=2\sigma \mathscr H^2(K)
\]
where $\sigma>0$ is the surface tension. 
\begin{figure}[htbp]
\begin{center}
\includegraphics[width=0.62\textwidth]{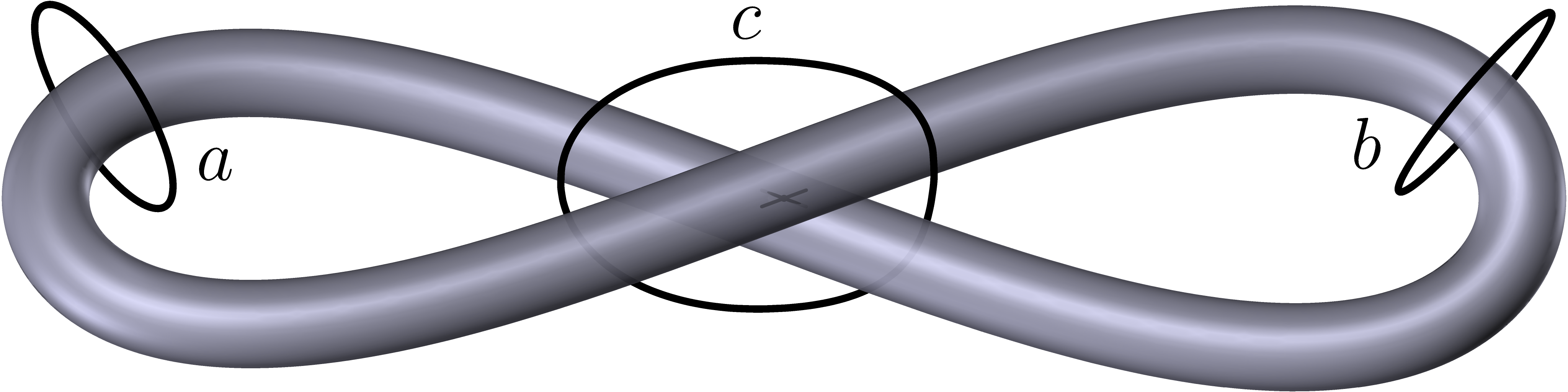}
\end{center}
\caption{For this loop, 
if we look for a spanning set relative to the homotopy class of the loops $a$ or $b$, spanning surfaces covering only the hole on the left or on the right one will be allowed, respectively. If, instead, we consider the homotopy class of the loop $c$, both holes must be covered by the spanning set.}
\label{fig:spanning}
\end{figure}

\begin{figure}[htbp]
		\centering
		\includegraphics[width=6cm]{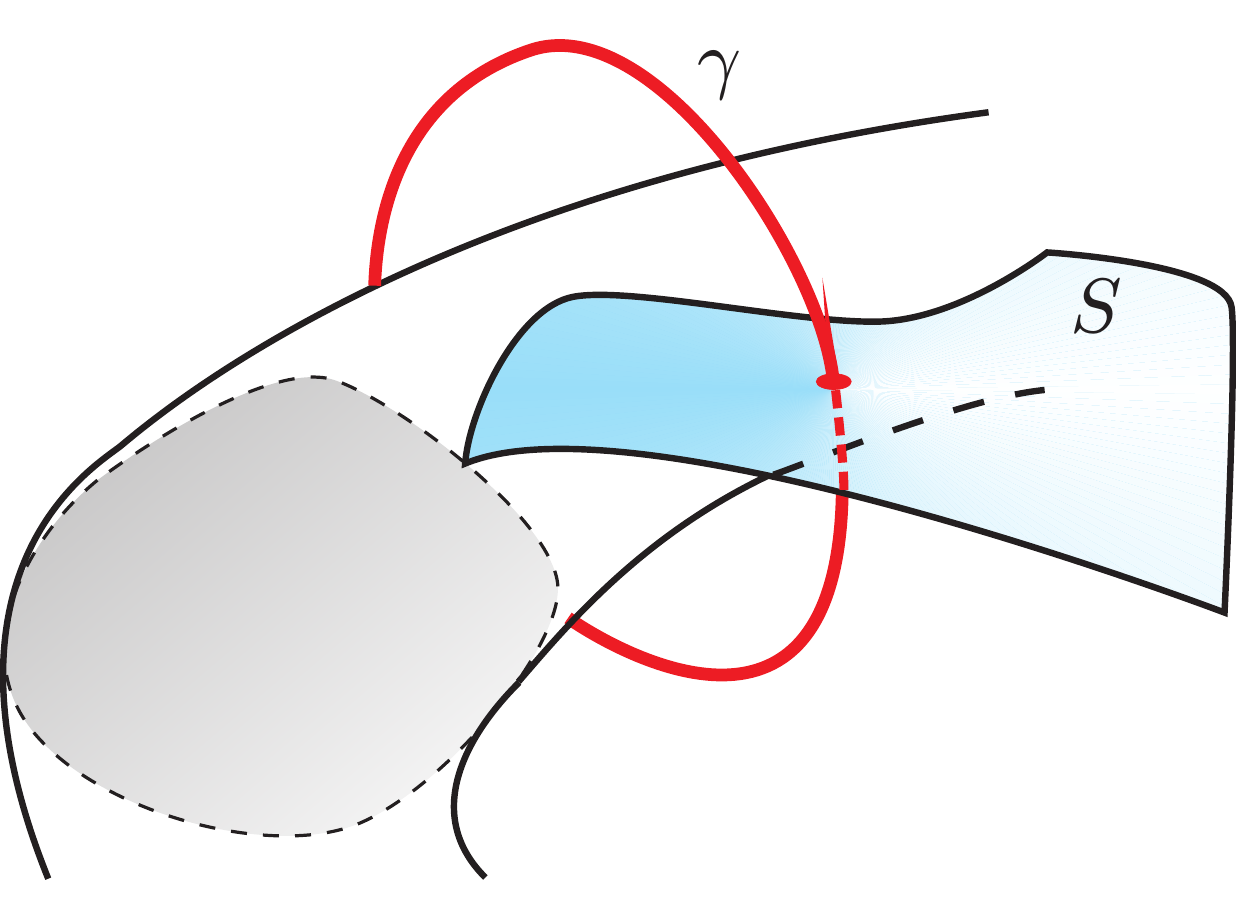}
	\caption{The surface $S$ must intersect the loop $\gamma$.}
	\label{fig:spanning_insieme}
\end{figure}

Using the framework and the notation presented in \cref{sec:minimalsets}, we impose the spanning condition choosing 
a suitable class of loops closed by homotopy. Indeed, an appropriate choice of homotopy classes determines which holes of a bounding loop with points of self-contact are covered, see \cref{fig:spanning}. Precisely, we use the subset $\mathscr D_{\Lambda[\vect w]} \subset \mathscr C_{\Lambda[\vect w]}$ containing all $\gamma$ that have linking number 1 or $-1$ with the midline $\vect x$. Then, we seek a surface $K\in \mathscr F(\Lambda[\vect w],\mathscr D_{\Lambda[\vect w]})$ that is a \mbox{$\mathscr D_{\Lambda[\vect w]}$-spanning} set of the bounding loop $\Lambda[\vect w]$ in the sense of \eqref{e:spanning} where $\mathcal C = \mathscr D_{\Lambda[\vect w]}$, see \cref{fig:spanning_insieme}.

The Kirchhoff-Plateau problem concerns the minimization of the energy functional 
\[
E_\mathrm{KP}(\vect w)=E_\mathrm{loop}(\vect w)+\inf\big\{E_\mathrm{film}(K) : K\in \mathscr F(\Lambda[\vect w],\mathscr D_{\Lambda[\vect w]})\big\}
\]
under the constraints \eqref{eq:closure}--\eqref{eq:glob-inj}. The main Theorem proved in \cite{GLF} is the following one.

\begin{theorem}\label{thm:existence}
Assume that there exists $\tilde{\vect w} \in V$ satisfying \eqref{eq:closure}--\eqref{eq:glob-inj} with $E_\mathrm{KP}(\tilde {\vect w})<+\infty$. Then there exists a minimizer of $E_\mathrm{KP}$ satisfying \eqref{eq:closure}--\eqref{eq:glob-inj}. Moreover, there is a relatively closed subset $K[\vect w]$ of $\R^3\setminus\Lambda[\vect w]$ such that 
\[
E_\mathrm{film}(K[\vect w])=\inf\big\{E_\mathrm{film}(K) : K\in \mathscr F(\Lambda[\vect w],\mathscr D_{\Lambda[\vect w]})\big\}.
\]
Finally, $K[\vect w]$ is a minimal set in $\R^3\setminus\Lambda[\vect w]$.
\end{theorem}

\noindent
{\it Idea of the proof.} Consider a minimizing sequence $(\vect w_h)$ for $E_{\rm KP}$ such that $E_{\rm KP}(\vect w_h) \le M$ for some $M\ge 0$. It is possible to extract a weakly converging subsequence, not relabeled, $\vect w_h\rightharpoonup  \vect w$ where $\vect w$ satisfies the constraints \eqref{eq:closure}--\eqref{eq:glob-inj}. 
The key point is to prove that if $K_h\in \mathscr F(\Lambda[\vect w_h],\mathcal D_{\Lambda[\vect w_h]})$ and a loop $\gamma$ in $\mathcal D_{\Lambda[\vect w]}$ then for any $\varepsilon>0$ such that the tubular neighborhood $U_{2\varepsilon}(\gamma)$ of radius $2\e$ around $\gamma$ is contained in $\mathbb R^3\setminus \Lambda[\vect w_h]$, there exists $M=M(\varepsilon)>0$ such that, for any $h$ large enough,
\begin{equation}\label{key}
\mathscr H^2(K_h \cap U_\varepsilon(\gamma))\ge M.
\end{equation}
Indeed, take $K_h$ with 
\[
\mathscr H^2(K_h)=\inf\big\{E_\mathrm{film}(K) : K_h\in \mathscr F(\Lambda[\vect w_h],\mathcal D_{\Lambda[\vect w_h]})\big\}.
\]
This is always possible essentially thanks to \cite[Thm.\,2]{DGM}. The measures $\mu_h:=\mathscr H^2\res K_h$ constitute a bounded sequence, $\mu_h\stackrel{*}{\rightharpoonup}\mu$ up to the extraction of a subsequence, and the limit measure satisfies
\[
\mu\geq\mathscr H^2\res K_\infty
\] 
where $K_{\infty}:={\rm spt}(\mu)\setminus\Lambda[\vect w]$ is a countably $\mathscr H^2$-rectifiable set. Assume by contradiction that there exists $\gamma \in \mathcal D_{\Lambda[\vect w]}$ with $\gamma \cap K_\infty= \emptyset$ and take $\varepsilon$ as before. We therefore find that $\mu(U_{2\varepsilon}(\gamma))=0$ and then 
\[
\lim_{h\to+\infty}\mathscr  H^2(K_h \cap U_{\varepsilon}(\gamma))=0
\]
which contradicts \eqref{key}. This means that $K_\infty\in \mathscr F(\Lambda[\vect w],\mathcal D_{\Lambda[\vect w]})$. 

We also get
\[
\begin{aligned}
\liminf_{h\to+\infty}&\inf\{\mathscr H^2(K) : K_h\in \mathscr F(\Lambda[\vect w_h],\mathcal D_{\Lambda[\vect w_h]})\}\\
&\ge \liminf_{h\to+\infty}\mathscr H^2(K_h)\\
&=\liminf_{h\to+\infty}\mathcal \mu_h(\mathbb R^3)\\
&\ge \mu(\mathbb R^3)\\
&\ge \mathscr H^2(K_\infty)\\
&\ge \inf\{\mathscr H^2(K) : K\in \mathscr F(\Lambda[\vect w],\mathcal D_{\Lambda[\vect w]})\}
\end{aligned}
\]
which establishes the lower semicontinuity of the functional $E_\mathrm{KP}$ (the lower semicontinuity of the loop energy is a standard).

Finally, the fact that there exists a minimal set $K[\vect w]$ in $\R^3\setminus\Lambda[\vect w]$  with 
\[
E_\mathrm{film}(K[\vect w])=\inf\big\{E_\mathrm{film}(K): K\in \mathscr F(\Lambda[\vect w],\mathcal D_{\Lambda[\vect w]})\big\}
\]
follows from \cite{DGM} and this yields the conclusion.
\qed

\subsection{Linked rods}

\begin{figure} [htbp]
\centering
\includegraphics[width=0.4\textwidth]{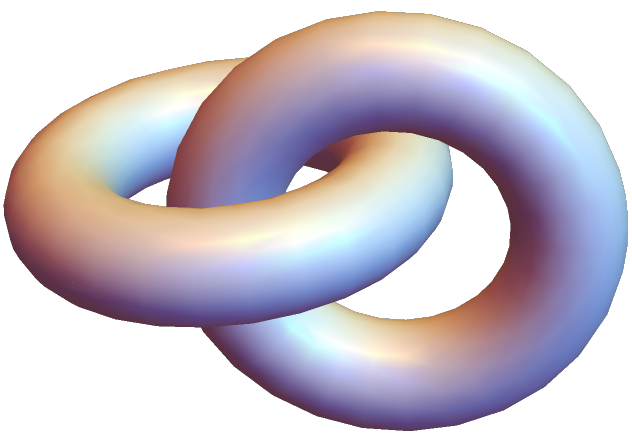}
\caption{A possible geometry of two linked rods.}
\label{sistema}
\end{figure}

The first generalization of the Kirchhoff-Plateau problem has been investigated in \cite{BLM1} where a more complex configuration of the bounding loop is considered. Precisely, the loop consists in a finite number of rods linked in an arbitrary way, as for instance in \cref{sistema}.
Following the same notation as before, and limiting to the case of two rods, two vectors $\vect w_1,\vect w_2$ are introduced to describe the two midlines: 
\[
\vect w^{(1)}=((\kappa_1^{(1)},\kappa_2^{(1)},\omega^{(1)}),\vect x_0^{(1)},\vect t_0^{(1)},\vect d_0^{(1)}), \quad \vect w^{(2)}=((\kappa_1^{(2)},\kappa_2^{(2)},\omega^{(2)}),\vect x_0^{(2)},\vect t_0^{(2)},\vect d_0^{(2)}).
\]
In particular, it is assumed that only the midline generated by $\vect w^{(1)}$ is clamped, that means that $(\vect x_0^{(1)},\vect t_0^{(1)},\vect d_0^{(1)})$ is prescribed. Concerning the second rod, we do not assume a priori its position in space, namely the vector $(\vect x_0^{(2)},\vect t_0^{(2)},\vect d_0^{(2)})$ is an unknown of the problem. 
\begin{figure}[htbp]
\centering
\includegraphics[width=0.68\textwidth]{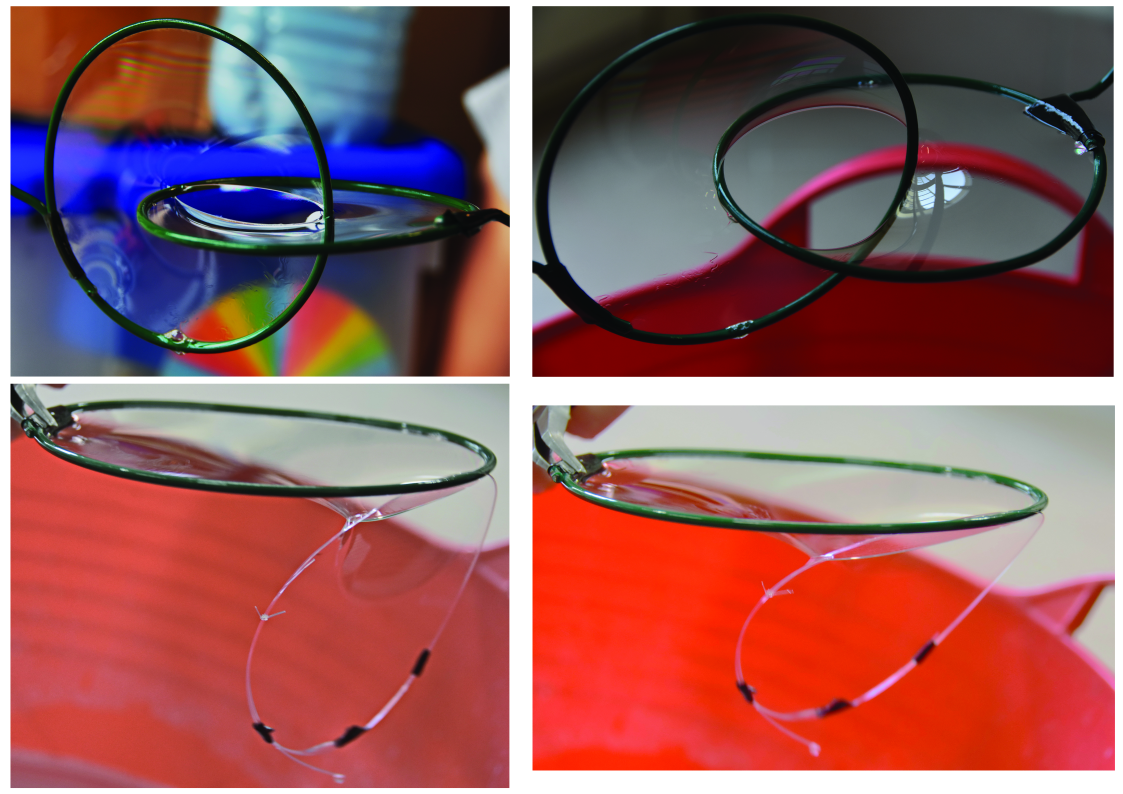}
\caption{Above: two fixed linked rigid metallic wires in a soap solution. Below: one rod is more flexible than the other one.}
\label{experiments1}
\end{figure}
For each rod we assume the corresponding analogous constraints \eqref{eq:closure}--\eqref{eq:glob-inj}. Moreover, we have also to ask that the linking number between the two midlines is prescribed:
\begin{equation}\label{link}
{\rm Link}(\vect x^{(1)},\vect x^{(2)})=\eta
\end{equation}
for some $\eta\in \Z$. Finally, since we need to have a non-interpenetration between the two rods we are going to assume that 
\begin{equation}\label{two}
\mathscr L^3(\Lambda[\vect w^{(1)}]\cap \Lambda[\vect w^{(2)}])=0.
\end{equation}
Concerning the spanning conditions, in this case we choose the loops which are not homotopic to a constant and such that the sum of the linking numbers between the loop and the two rods is always one: this means that a loop cannot link at the same time both rods. Thus, the Kirchhoff-Plateau problem concerns the minimization of the following energy functional
\[
(\vect w^{(1)},\vect w^{(2)}) \mapsto E_\mathrm{loop}(\vect w^{(1)})+E_\mathrm{loop}(\vect w^{(2)})+\inf\big\{E_\mathrm{film}(K) : \textrm{$K$ spans $\Lambda[\vect w^{(1)}]\cup \Lambda[\vect w^{(2)}]$})\big\},
\]
under all of the constraints described above. In \cite{BLM1}, we provide the existence of a minimizer and we perform some experiments, see \cref{experiments1}.


\subsection{Soap films spanning repulsive links}

\begin{figure}[htbp]
\centering
\includegraphics[width=0.7\textwidth]{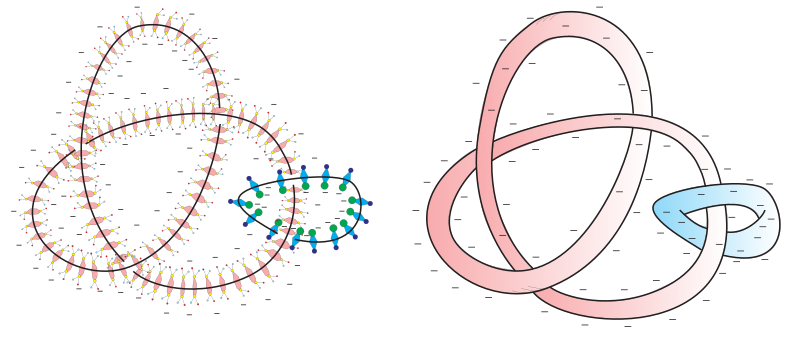}
\caption{Knotted protein linked to another one.}
\label{biosistema}
\end{figure}

As a second generalization, in \cite{BLM2} the case of knotted proteins is treated. To consider processes like the adsorption of a protein by a biomembrane, in \cite{BLM2} we introduce an additional repulsional energy between the two linked rods; see \cref{biosistema}.

\begin{figure}[htbp]
\centering
\includegraphics[width=0.55\textwidth]{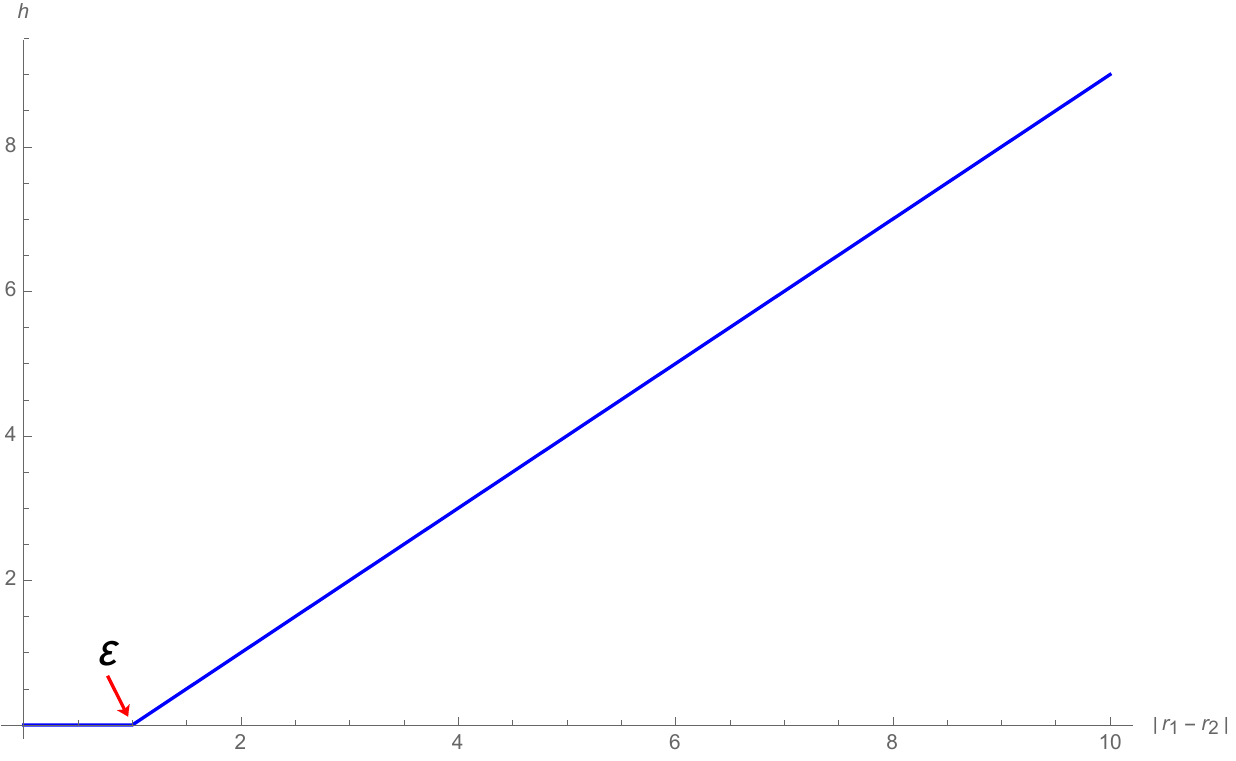}
\caption{Example of the ``repulsive'' function $h$.}
\label{fig:h}
\end{figure}

The general setting is the same as in the case previously considered of two linked rods: $\vect w^{(1)}$ and $\vect w^{(2)}$ generate two midlines $\vect x^{(1)}$ and $\vect x^{(1)}$ respectively. We assume all the usual constraints \eqref{eq:closure}--\eqref{eq:glob-inj} on $\vect x^{(i)}$, as well as \eqref{link}. We substitute \eqref{two} with electrical potential energy term which, physically, encodes the repulsion between the two rods. Precisely, the repulsion is modeled by  
\begin{equation}
\label{erep}
\int_{0}^{L_1}\int_{0}^{L_2} \frac{1}{h(\|\vect x^{(1)}(s_1)- \vect x^{(2)}(s_2)\|)}\, ds_1ds_2,
\end{equation}
where $h$ is a suitable increasing, nonnegative and continuous function. With this choice, we are introducing a positively unbounded energy, that may be infinite if the midlines are sufficiently close. 
A possible choice for $h$ is represented in \cref{fig:h}: a function which is $0$ until some positive and small parameter $\e$ and then grows linearly. Therefore, the energy functional becomes 
\[
\begin{aligned}
(\vect w^{(1)},\vect w^{(2)}) \mapsto E_\mathrm{loop}(\vect w^{(1)})+E_\mathrm{loop}(\vect w^{(2)})&+\int_{0}^{L_1}\int_{0}^{L_2} \frac{1}{h(\|\vect x^{(1)}(s_1)- \vect x^{(2)}(s_2)\|)}\, ds_1ds_2\\
&+\inf\big\{E_\mathrm{film}(K) : \textrm{$K$ spans $\Lambda[\vect w^{(1)}]\cup \Lambda[\vect w^{(2)}]$})\big\}
\end{aligned}
\]
which is minimized under all of the described constraints.

\subsection{Dimensional reduction}
The paper \cite{BLM3} deals with a sort of dimensional reduction of the Kirchhoff-Plateau problem. The aim of that work is then to perform a formal dimensional reduction of the classical Kirchhoff-Plateau problem where the limiting curve is the midline of the rod. We require that $\vect w$ satisfies the assumptions \eqref{eq:closure}--\eqref{eq:glob-inj}. 
Moreover, we assume that the so-called {\it global radius of curvature} $\Delta(\vect x[\vect w])]$ is bounded from below by a constant $\Delta_0>0$; this assumption prevents self-intersection for a sufficiently small cross section (see \cite{gms}). In this context the cross section $A(s)$ is replaced by a rescaled cross section $\e A(s)$, where $\e>0$ is a positive and vanishing parameter. In addition, we can define the map $\vect p^\e[\vect w]$ and the corresponding rod $\Lambda^\e[\vect w]=\vect p^\e[\vect w](\Omega^\e)$. The rescaled energy functional reads as 
\[
E_\mathrm{KP}^\e(\vect w)=E_\mathrm{loop}^\e(\vect w)+\inf\big\{E_\mathrm{film}(K) : K\in \mathscr F(\Lambda^\e[\vect w],\mathscr D_{\Lambda^\e[\vect w]})\big\}
\]
where 
\[
E_\mathrm{loop}^\e(\vect w)=E_{\rm sh}(\vect w)-\frac{1}{\e^2}\int_{\Omega^\e}\rho(s,\zeta_1,\zeta_2)\vect g\cdot \vect p^\e[\vect w](s,\zeta_1,\zeta_2)\,dsd\zeta_1d\zeta_2.
\]
As $\e\to 0^+$ we obtain a limit functional (in the sense of $\Gamma$-convergence) which is given by 
\[
E^0(\vect w)=E_{\rm sh}(\vect w)-\int_0^L|A(s)|\rho_0(s)\vect g\cdot \vect x[\vect w](s)\,ds+\inf\big\{E_\mathrm{film}(K) : \textrm{$K$ spans $\vect x[\vect w]([0,L])$}\big\}
\]
being 
\[
\rho_0(s)=\lim_{(\xi_1,\xi_2)\to (0,0)}\rho(s,\xi_1,\xi_2).
\]
The approximating problems have minima which converge weakly to the minimum energy solution of the limit problem, as well as the corresponding value of the energy. This also shows that the Plateau solution with elastic line boundary may be approximated by solutions of the problems with a rod boundary.

\section{Further results and work in progress}\label{sec-working_progress}

Once the existence of minimizers has been obtained, it is useful to characterize them deriving, for instance, the Euler-Lagrange equations. Unfortunately, since the definition of the bounding loop requires a high number of constraints, a first simplification is to consider an elastic curve instead of the Kirchhoff-rod as the boundary wire. Thus, a first simplification can be considering the Elastic-Plateau problem: we are interested in performing the variational analysis of energy functionals of the type
\begin{equation*}
\mathcal{E}[\vect{\gamma}, \vect{X}] = \int_{\vect{\gamma}} f(\kappa, \tau)\, d\ell + {\rm Area}(S)
\end{equation*}
where $\vect{\gamma}$ is a closed curve in $\R^3$ with curvature $\kappa$ and torsion $\tau$ and ${\rm Area}(S)$ is the area of the spanning surface spanning the elastic curve $\vect{\gamma}$. The derivation of minimizers and their characterization leads to several difficulties like getting compactness in the disc-type approach or dealing with Plateau singularities. It seems that one should use the framework of Lytchak and Wenger \cite{lytchak2017area} and Creutz \cite{creutz2022plateau} to set the problem in Sobolev spaces for disc-type surfaces spanning a curve with possible self-intersections, while we expect to deal with Geometric Measure Theory to treat general surfaces.

A first attempt to investigate the mentioned problem has been done in \cite{BLM_PRSA, BLM_INDAM, BBLM2024}, where in order to deal with elastic curves, we minimize $\mathcal{E}$ among all disc-type maps $\vect X\colon D\to \R^3$ with trace pointwise equal to the elastic curve $\vect{\gamma}$ (this condition is different from the classical Plateau problem rigorously solved by Douglas and Rad\`{o} \cite{Douglas, Rado} where the trace of $\vect{X}$ is a suitable reparametrization $\sigma(s)$ with $s \in [0,1]$ of the curve $\gamma$). Moreover, the area functional is substituted with $\int_D \Psi(\nabla \vect{X})\, dudv$,
where $D$ is the unit disc in $\R^2$ and $\vect X\colon D \to \R^3$ is a parametrization of a membrane spanning the elastic curve $\vect\gamma$, modeled through its deformation gradient $\nabla \vect{X} \in \R^{3\times2}$.
We adopt two different approaches to model the line integral: {\em the parametrized curves approach and the framed curves approach}. For the first one, the curve $\vect \gamma$ is modeled as the Euler-Bernoulli elastica, while only linear elastic membranes are taken into account (For details we refer to \cite[Theorem 2.2 - Theorem 2.3 - Theorem 2.7]{BBLM2024}). Concerning the second one, it is introduced to deal with more general energies, both for the boundary curve and for the membrane. Precisely, we introduce a moving orthonormal frame $\left\{\vect{t},\vect{n},\vect{b}\right\} \in W^{1, p}((0, 2 \pi); SO(3))$ with $p >1$ which generates a curve $\vect{r}$ by integration. On this basis, we impose suitable constraints to get a closed curve (For details we refer to \cite[Theorem 3.1]{BLM_PRSA} and \cite[Theorem 3.3 - Theorem 3.5]{BBLM2024}). 

Moreover, another interesting direction of investigation would be to perform a numerical study in order to visualize minimizers and their behaviour. A first attempt has been proposed in \cite[Section 4]{BBLM2024} where {\em ad hoc} method has been introduced to test some simple configurations in the membrane case. Precisely, developing a numerical approach is quite challenging due to the large number of constraints in the formulation of the problem, for instance we mention the pointwise length preserving constraint which is ill-suited in the application of a finite element method, or the choice of energy functionals, non linear terms are hard to be treat numerically. 

Finally, in its classical formulation, the Plateau problem is an {optimization problem}: looking for the surface with minimal area spanning the assigned boundary. However, it would be interesting to characterize the dynamical process since, especially from the physical viewpoint, Plateau devised many experiments putting a wire frame into a soap solution to a soap film.  
In particular, a first step can be to formulate and solve the dynamical Plateau problem in its quasi-static approximation: the idea is to prescribe the motion of the elastic curve and, at each time step $t \in [0,1]$, a minimal surface spanning the assigned curve must be determined. This approach generalizes the machinery introduced by Dal Maso and co-authors for studying fractures \cite{dal2002model}.

\section*{Acknowledgement \& Funding}
GB is supported by the European Research Council (ERC), under the European Union's Horizon 2020 research and innovation program, through the project ERC VAREG - {\em Variational approach to the regularity of the free boundaries} (grant agreement No. 853404) and GB acknowledges the MIUR Excellence Department Project awarded to the Department of Mathematics, University of Pisa, CUP I57G22000700001. GB and LL are supported by Gruppo Nazionale per l'Analisi Matematica, la Probabilit\`a e le loro Applicazioni (GNAMPA) of Istituto Nazionale di Alta Matematica (INdAM) though the INdAM-GNAMPA project 2024 CUP E53C23001670001. 
AM is supported by Gruppo Nazionale per la Fisica Matematica (GNFM) of Istituto Nazionale di Alta Matematica (INdAM).

\printbibliography

@book {AFP,
    AUTHOR = {Ambrosio, Luigi and Fusco, Nicola and Pallara, Diego},
     TITLE = {Functions of bounded variation and free discontinuity
              problems},
    SERIES = {Oxford Mathematical Monographs},
 PUBLISHER = {The Clarendon Press, Oxford University Press, New York},
      YEAR = {2000},
     PAGES = {xviii+434},
}

@article{creutz2022plateau,
  title={Plateau’s problem for singular curves},
  author={Creutz, Paul},
  journal={Communications in Analysis and Geometry},
  volume={30},
  number={8},
  pages={1779--1792},
  year={2022},
  publisher={International Press of Boston}
}

@article{lytchak2017area,
  title={Area minimizing discs in metric spaces},
  author={Lytchak, Alexander and Wenger, Stefan},
  journal={Archive for Rational Mechanics and Analysis},
  volume={223},
  pages={1123--1182},
  year={2017},
  publisher={Springer}
}

@article{HoaFri16,
  title={Influence of a spanning liquid film on the stability and buckling of a circular loop with intrinsic curvature or intrinsic twist density},
  author={Hoang, Tuan and Fried, Eliot},
  journal={Mathematics and Mechanics of Solids},
  volume={23},
  number={1},
  pages={43--66},
  year={2018},
  publisher={SAGE Publications Sage UK: London, England}
}

@article{bevilacqua2024effects,
  title={Effects of surface tension and elasticity on critical points of the Kirchhoff--Plateau problem},
  author={Bevilacqua, Giulia and Lonati, Chiara},
  journal={Bollettino dell'Unione Matematica Italiana},
  volume={17},
  number={2},
  pages={221--240},
  year={2024},
  publisher={Springer}
}

@article{CheFri14,
  title={Stability and bifurcation of a soap film spanning a flexible loop},
  author={Chen, Yi-chao and Fried, Eliot},
  journal={Journal of Elasticity},
  volume={116},
  pages={75--100},
  year={2014},
  publisher={Springer}
}

@article{BirFri15,
  title={Theoretical and experimental study of the stability of a soap film spanning a flexible loop},
  author={Biria, Aisa and Fried, Eliot},
  journal={International Journal of Engineering Science},
  volume={94},
  pages={86--102},
  year={2015},
  publisher={Elsevier}
}

@article{BirFri14,
  title={Buckling of a soap film spanning a flexible loop resistant to bending and twisting},
  author={Biria, Aisa and Fried, Eliot},
  journal={Proceedings of the Royal Society A: Mathematical, Physical and Engineering Sciences},
  volume={470},
  number={2172},
  pages={20140368},
  year={2014},
  publisher={The Royal Society Publishing}
}

@article{GFF,
  title={Instability paths in the Kirchhoff--Plateau problem},
  author={Giusteri, Giulio G and Franceschini, Paolo and Fried, Eliot},
  journal={Journal of Nonlinear Science},
  volume={26},
  pages={1097--1132},
  year={2016},
  publisher={Springer}
}

@article{BY,
  title={Minimal surfaces with an elastic boundary},
  author={Bernatzki, Felicia and Ye, Rugang},
  journal={Annals of Global Analysis and Geometry},
  volume={19},
  number={1},
  pages={1--9},
  year={2001},
  publisher={Kluwer Academic Publishers Dordrecht}
}

@article{B,
  title={On the existence and regularity of mass-minimizing currents with an elastic boundary},
  author={Bernatzki, Felicia},
  journal={Annals of Global Analysis and Geometry},
  volume={15},
  pages={379--399},
  year={1997},
  publisher={Springer}
}

@book{do2016differential,
  title={Differential geometry of curves and surfaces: revised and updated second edition},
  author={Do Carmo, Manfredo P},
  year={2016},
  publisher={Courier Dover Publications}
}

@article{Rado,
  title={The problem of the least area and the problem of {P}lateau},
  author={Rad{\'o}, Tibor},
  journal={Mathematische Zeitschrift},
  volume={32},
  number={1},
  pages={763--796},
  year={1930},
  publisher={Springer}
}

@article{Douglas,
  title={Solution of the problem of {P}lateau},
  author={Douglas, Jesse},
  journal={Transactions of the American Mathematical Society},
  volume={33},
  number={1},
  pages={263--321},
  year={1931},
  publisher={JSTOR}
}

@book {D,
    AUTHOR = {Dacorogna, Bernard},
     TITLE = {Direct methods in the calculus of variations},
    SERIES = {Applied Mathematical Sciences},
    VOLUME = {78},
   EDITION = {Second},
 PUBLISHER = {Springer, New York},
      YEAR = {2008},
     PAGES = {xii+619},
}

@book{DHS,
  title={Minimal surfaces},
  author={Dierkes, Ulrich and Hildebrandt, Stefan and Sauvigny, Friedrich},
  year={2010},
  publisher={Springer}
}

@article{G,
  title={Regularity of minimizing surfaces of prescribed mean curvature},
  author={Gulliver, Robert D},
  journal={Annals of Mathematics},
  volume={97},
  number={2},
  pages={275--305},
  year={1973},
  publisher={JSTOR}
}

@article{dal2002model,
  title={A Model for the Quasi-Static Growth of Brittle Fractures: Existence and Approximation Results},
  author={Dal Maso, Gianni and Toader, Rodica},
  journal={Archive for Rational Mechanics and Analysis},
  volume={162},
  pages={101--135},
  year={2002},
  publisher={Springer}
}

@article{HS,
  title={On some non-linear elliptic differential-functional equations},
  author={Hartman, Philip and Stampacchia, Guido},
  journal={Acta Mathematica},
  volume={115},
  number={1},
  pages={271--310},
  year={1966},
  publisher={Springer}
}

@book{K,
  title={{\"U}ber Minimalfl{\"a}chen, deren Randkurven wenig von ebenen Kurven abweichen},
  author={Korn, Arthur},
 journal={Phys.-Math. Cl.},
volume={II},
  pages={1--37},
  year={1909},
  publisher={Verlag der K{\"o}nigl. Akad. der Wiss.}
}

@article{JS,
  title={Variational problems of minimal surface type II. Boundary value problems for the minimal surface equation},
  author={Jenkins, Howard and Serrin, James},
  journal={Archive for Rational Mechanics and Analysis},
  volume={21},
  number={4},
  pages={321--342},
  year={1966},
  publisher={Springer}
}

@book{Maggi,
  title={Sets of finite perimeter and geometric variational problems: an introduction to Geometric Measure Theory},
  author={Maggi, Francesco},
  number={135},
  year={2012},
  publisher={Cambridge University Press}
}

@article{Mo,
  title={On the solutions of quasi-linear elliptic partial differential equations},
  author={Morrey, Charles B},
  journal={Transactions of the American Mathematical Society},
  volume={43},
  number={1},
  pages={126--166},
  year={1938},
  publisher={JSTOR}
}

@article{Mu,
  title={Zum Randwertproblem der partiellen Differentialgleichung der Minimalfl{\"a}chen.},
  author={M{\"u}ntz, Ch},
  journal={Journal f\"{u}r die {R}eine und {A}ngewandte {M}athematik},
  volume={139},
 pages={52--79},
  year={1911},
  publisher={Walter de Gruyter, Berlin/New York Berlin, New York}
}

@book{N,
  title={Lecture on minimal surfaces},
  author={Nitsche, Johannes C},
  year={1989},
  publisher={Cambridge university press}
}

@article{SS,
  title={Derivations of the Young-Laplace equation},
  author={Siqveland, Leiv Magne and Skj{\ae}veland, Svein Magne},
  journal={Capillarity},
  volume={4},
  number={2},
  pages={23--30},
  year={2021}
}

@book {S,
    AUTHOR = {Simon, Leon},
     TITLE = {Lectures on geometric measure theory},
    SERIES = {Proceedings of the Centre for Mathematical Analysis,
              Australian National University},
    VOLUME = {3},
 PUBLISHER = {Australian National University, Centre for Mathematical
              Analysis, Canberra},
      YEAR = {1983},
     PAGES = {vii+272},
      ISBN = {0-86784-429-9},
   MRCLASS = {49-01 (28A75 49F20)},
  MRNUMBER = {756417},
MRREVIEWER = {J.\ S.\ Joel},
}

@book {F,
    AUTHOR = {Federer, Herbert},
     TITLE = {Geometric measure theory},
    SERIES = {Die Grundlehren der mathematischen Wissenschaften},
    VOLUME = {Band 153},
 PUBLISHER = {Springer-Verlag New York, Inc., New York},
      YEAR = {1969},
     PAGES = {xiv+676},
   MRCLASS = {28.80 (26.00)},
  MRNUMBER = {257325},
MRREVIEWER = {J.\ E.\ Brothers},
}

@book {A,
    AUTHOR = {Almgren, F.J.},
     TITLE = {The theory of varifolds},
    SERIES = {},
    VOLUME = {},
 PUBLISHER = {Mimeographed notes, Princeton},
      YEAR = {1965},
     PAGES = {},
   MRCLASS = {},
  MRNUMBER = {},
MRREVIEWER = {},
}

@article {H,
    AUTHOR = {Hutchinson, John E.},
     TITLE = {Second fundamental form for varifolds and the existence of
              surfaces minimising curvature},
   JOURNAL = {Indiana Univ. Math. J.},
  FJOURNAL = {Indiana University Mathematics Journal},
    VOLUME = {35},
      YEAR = {1986},
    NUMBER = {1},
     PAGES = {45--71},
}

@article {bdgg,
    AUTHOR = {Bombieri, E. and De Giorgi, E. and Giusti, E.},
     TITLE = {Minimal cones and the {B}ernstein problem},
   JOURNAL = {Invent. Math.},
  FJOURNAL = {Inventiones Mathematicae},
    VOLUME = {7},
      YEAR = {1969},
     PAGES = {243--268},
}

@article {A1,
    AUTHOR = {Almgren, Jr., F. J.},
     TITLE = {Existence and regularity almost everywhere of solutions to
              elliptic variational problems with constraints},
   JOURNAL = {Mem. Amer. Math. Soc.},
  FJOURNAL = {Memoirs of the American Mathematical Society},
    VOLUME = {4},
      YEAR = {1976},
    NUMBER = {165},
     PAGES = {viii+199},
}

@article {T,
    AUTHOR = {Taylor, Jean E.},
     TITLE = {The structure of singularities in soap-bubble-like and
              soap-film-like minimal surfaces},
   JOURNAL = {Ann. of Math.},
  FJOURNAL = {Annals of Mathematics. Second Series},
    VOLUME = {103},
      YEAR = {1976},
    NUMBER = {3},
     PAGES = {489--539},
}

@article {HP,
    AUTHOR = {Harrison, Jenny and Pugh, Harrison},
     TITLE = {Existence and soap film regularity of solutions to {P}lateau's
              problem},
   JOURNAL = {Adv. Calc. Var.},
  FJOURNAL = {Advances in Calculus of Variations},
    VOLUME = {9},
      YEAR = {2016},
    NUMBER = {4},
     PAGES = {357--394},
}

@article {Ha,
    AUTHOR = {Harrison, J.},
     TITLE = {Soap film solutions to {P}lateau's problem},
   JOURNAL = {J. Geom. Anal.},
  FJOURNAL = {Journal of Geometric Analysis},
    VOLUME = {24},
      YEAR = {2014},
    NUMBER = {1},
     PAGES = {271--297},
}

@article {DGM,
    AUTHOR = {De Lellis, C. and Ghiraldin, F. and Maggi, F.},
     TITLE = {A direct approach to {P}lateau's problem},
   JOURNAL = {J. Eur. Math. Soc. (JEMS)},
  FJOURNAL = {Journal of the European Mathematical Society},
    VOLUME = {19},
      YEAR = {2017},
    NUMBER = {8},
     PAGES = {2219--2240},
}

@article {GLF,
    AUTHOR = {Giusteri, Giulio G. and Lussardi, Luca and Fried, Eliot},
     TITLE = {Solution of the {K}irchhoff-{P}lateau problem},
   JOURNAL = {J. Nonlinear Sci.},
  FJOURNAL = {Journal of Nonlinear Science},
    VOLUME = {27},
      YEAR = {2017},
    NUMBER = {3},
     PAGES = {1043--1063},
}

@article {BLM2,
    AUTHOR = {Bevilacqua, Giulia and Lussardi, Luca and Marzocchi, Alfredo},
     TITLE = {Soap film spanning electrically repulsive elastic protein
              links},
   JOURNAL = {Atti Accad. Peloritana Pericolanti Cl. Sci. Fis. Mat. Natur.},
  FJOURNAL = {Atti della Accademia Peloritana dei Pericolanti. Classe di
              Scienze, Fisiche, Matematiche e Naturali. AAPP. Physical,
              Mathematical, and Natural Sciences},
    VOLUME = {96},
      YEAR = {2018},
     PAGES = {A1, 13},
}

@article {BLM1,
    AUTHOR = {Bevilacqua, Giulia and Lussardi, Luca and Marzocchi, Alfredo},
     TITLE = {Soap film spanning an elastic link},
   JOURNAL = {Quart. Appl. Math.},
  FJOURNAL = {Quarterly of Applied Mathematics},
    VOLUME = {77},
      YEAR = {2019},
    NUMBER = {3},
     PAGES = {507--523},
}

@article {BLM3,
    AUTHOR = {Bevilacqua, Giulia and Lussardi, Luca and Marzocchi, Alfredo},
     TITLE = {Dimensional reduction of the {K}irchhoff-{P}lateau problem},
   JOURNAL = {J. Elasticity},
  FJOURNAL = {Journal of Elasticity. The Physical and Mathematical Science
              of Solids},
    VOLUME = {140},
      YEAR = {2020},
    NUMBER = {1},
     PAGES = {135--148},
}

@article{BBLM2024,
  title={Elastic membranes spanning deformable curves},
  author={Ballarin, Francesco and Bevilacqua, Giulia and Lussardi, Luca and Marzocchi, Alfredo},
  journal={ZAMM-Journal of Applied Mathematics and Mechanics/Zeitschrift f{\"u}r Angewandte Mathematik und Mechanik},
  pages={e202300890},
YEAR = {2024},
  publisher={Wiley Online Library},
}

@article{BLM_INDAM,
  title={Geometric invariants of non-smooth framed curves},
  author={Bevilacqua, Giulia and Lussardi, Luca and Marzocchi, Alfredo},
  journal={Accepted in INdAM-Series, ``Anisotropic Isoperimetric Problems \& Related Topics", arXiv preprint arXiv:2301.03525},
  year={2023}
}

@article{BLM_PRSA,
  title={Variational analysis of inextensible elastic curves},
author={Bevilacqua, Giulia and Lussardi, Luca and Marzocchi, Alfredo},
  journal={Proceedings of the Royal Society A},
  volume={478},
  number={2260},
  pages={20210741},
  year={2022},
  publisher={The Royal Society}
}

@article {DL,
    AUTHOR = {De Rosa, Antonio and Lussardi, Luca},
     TITLE = {On the anisotropic {K}irchhoff-{P}lateau problem},
   JOURNAL = {Math. Eng.},
  FJOURNAL = {Mathematics in Engineering},
    VOLUME = {4},
      YEAR = {2022},
    NUMBER = {2},
     PAGES = {Paper No. 011, 13},
}

@book {Antman2005,
    AUTHOR = {Antman, Stuart S.},
     TITLE = {Nonlinear problems of elasticity},
    SERIES = {Applied Mathematical Sciences},
    VOLUME = {107},
   EDITION = {Second},
 PUBLISHER = {Springer, New York},
      YEAR = {2005},
     PAGES = {xviii+831},
}

@book {Hartman1982,
    AUTHOR = {Hartman, Philip},
     TITLE = {Ordinary differential equations},
   EDITION = {second},
 PUBLISHER = {Birkh\"auser, Boston, MA},
      YEAR = {1982},
     PAGES = {xv+612},
}

@article {gms,
    AUTHOR = {Gonzalez, O. and Maddocks, J. H. and Schuricht, F. and von der
              Mosel, H.},
     TITLE = {Global curvature and self-contact of nonlinearly elastic
              curves and rods},
   JOURNAL = {Calc. Var. Partial Differential Equations},
  FJOURNAL = {Calculus of Variations and Partial Differential Equations},
    VOLUME = {14},
      YEAR = {2002},
    NUMBER = {1},
     PAGES = {29--68},
}

\end{document}